\documentclass[1p]{elsarticle}

\usepackage[english]{babel}

\usepackage{amsmath}
\usepackage{amssymb}
\usepackage{bm}
\usepackage{mathtools}
\usepackage{gensymb}

\usepackage{graphicx}
\graphicspath{{figures/}}

\usepackage{subcaption}
\usepackage{xspace}

\usepackage{nameref}
\makeatletter
\newcommand*{\currentname}{\@currentlabelname}
\makeatother

\DeclareMathOperator*{\argmax}{arg\,max}

\newcommand{\bigO}{\mathcal{O}}
\newcommand{\norm}[2]{\left\| #1 \right\|_{#2}}
\renewcommand{\vec}[1]{\ensuremath{\bm{#1}}}
\newcommand{\pref}[1]{(\ref{#1})}

\usepackage{booktabs}

\renewcommand{\epsilon}{\varepsilon}

\usepackage[section]{algorithm}
\usepackage{algorithmic}

\usepackage{setspace}

\usepackage{multirow}

\usepackage{lineno,hyperref}

\usepackage{xcolor}

\usepackage[section]{placeins}

\renewcommand{\cite}[1]{\citeauthor{#1} \citep{#1}}

\definecolor{revisedText}{rgb}{0,0,0.8}
\definecolor{revisedText}{rgb}{0,0,0}

\begin{document}

\begin{frontmatter}

\title{Implicit reconstructions of thin leaf surfaces from large, noisy point clouds}

\author[QUTaddress]{Riley M. Whebell}
\author[QUTaddress]{Timothy J. Moroney\corref{corauth}}
\ead{t.moroney@qut.edu.au}
\author[QUTaddress]{Ian W. Turner}
\author[QUTaddress]{Ravindra Pethiyagoda}
\author[QUTaddress]{Scott W. McCue}

\address[QUTaddress]{School of Mathematical Sciences, Queensland University of Technology, Brisbane QLD 4001, Australia}
\cortext[corauth]{Corresponding author}

\begin{abstract}
Thin surfaces, such as the leaves of a plant, pose a significant challenge for implicit surface reconstruction techniques, which typically assume a closed, orientable surface. 
We show that by approximately interpolating a point cloud of the surface (augmented with off-surface points) and restricting the evaluation of the interpolant to a tight domain around the point cloud, we need only require an orientable surface for the reconstruction.
We use polyharmonic smoothing splines to fit approximate interpolants to noisy data, and a partition of unity method with an octree-like strategy for choosing subdomains.
This method enables us to interpolate an $N$-point dataset in $\mathcal{O}(N)$ operations.
We present results for point clouds of capsicum and tomato plants, scanned with a handheld device.
An important outcome of the work is that sufficiently smooth leaf surfaces are generated that are amenable for droplet spreading simulations. 
\end{abstract}


\end{frontmatter}


\section{Introduction}

Surface reconstruction is the problem of constructing a digital representation of a physical object from scanned data. 
Techniques for reconstructing a surface with global smoothness are of specific interest in our related work on simulating droplet retention on plant surfaces \citep{mayo2015simulating,Dorr2014,Dorr2016,Zabkiewicz2020,Oqielat2011}.
Many algorithms have been proposed for this type of reconstruction from point cloud datasets \citep{Hoppe1992,Carr2001,Ohtake2003,Tobor2004,Xiaojun2005,Kazhdan2006,Zhong2019}.
\textcolor{revisedText}{
These algorithms typically consider watertight surfaces; 2D manifolds without boundary. 
The resolution required to reconstruct a plant leaf as a volume (with thickness) is prohibitively fine, so instead we imagine it to be infinitesimally thin; a 2D manifold with boundary.
}
This constraint, together with the often complex geometry of whole plants, poses a challenge for established surface reconstruction algorithms.
We will demonstrate that by casting the surface reconstruction problem as an approximate scattered data interpolation problem, a novel combination of well-established techniques can reconstruct globally smooth thin surfaces from large, noisy point clouds.

The datasets from which we reconstruct surfaces are point clouds generated by a handheld Artec Eva scanner \citep{arteceva}.
The scanner captures 16 frames of range data per second, which are then registered together by the included proprietary software.
An example of a scanned point cloud is shown in Figure \ref{fig:noisyPC}.
This data is quite noisy -- this is not entirely due to the measurement error of the Artec Eva scanner (accurate to 0.1mm), but also due to some slight movement of the plant while scanning, and the varying thickness of its leaves.

Many models for droplet motion on plant surfaces involve the surface curvature \citep{mayo2015simulating}, so our surface reconstructions should be at least twice continuously differentiable ($C^2$). Furthermore, we would like to distinguish between distinct leaves of the plant -- to study them individually or extract statistics to describe the plant -- so our reconstruction should be easily segmentable. 
Lastly, our method should be robust to typical large, noisy datasets generated by a handheld range scanner, so that they may be applied to scans collected by an end-user of spray retention modelling software (such as that developed by Dorr et al. \citep{Dorr2014,Dorr2016}).

\cite{Kempthorne2014} compared some techniques for fitting continuously differentiable ($C^1$) surfaces to scanned point clouds of various plant leaves.
\textcolor{revisedText}{They compared three explicit interpolation methods (discrete smoothing $D^2$ splines, the thin plate spline finite element smoother, and the radial basis function Clough-Tocher method) to fit surfaces of the form $z = \mathcal{F}(x,y)$; explicit interpolants.
While those methods resulted in smooth surfaces and reproduced the leaf morphology well}, there is a limitation to the explicit approach when applied to leaves with more complex structure.
For example, the surface may be multiply-defined in the $z$-dimension, but $\mathcal{F}(x,y)$ may only take on one value. 
\cite{Kempthorne2015} provided a solution to this problem \textcolor{revisedText}{for the long, thin leaves of wheat. They defined a new locally orthogonal coordinate system with respect to the leaf: $(u,v,w)$.
The $u$ direction was the local principal direction of the leaf surface, the $v$ direction was orthogonal to $u$ but still tangent to the leaf, and the $w$ direction was normal to the leaf surface.}
An explicit interpolant was then formed: $w = \mathcal{F}(u,v)$.
This approach, however, requires the construction of the coordinate system for each individual leaf, and manual segmentation of leaves from the plant (for individual fitting).
Other algorithms in computer graphics contexts have presented a solution to these parameterisation limitations using implicit surface reconstruction.

Implicit methods \citep{Hoppe1992,Carr2001,Ohtake2003,Kazhdan2006,Zhong2019} construct a scalar field $\mathcal{F}: \Omega \to \mathbb{R}$, where $\Omega \subseteq \mathbb{R}^3$ is an appropriate domain about the point cloud, such that the surface $\sigma$ is a level set of $\mathcal{F}$: 
\[ \sigma = \{ \vec{x}\in\Omega : \mathcal{F}(\vec{x}) = 0 \}. \]
These methods vary in their approach to constructing the function $\mathcal{F}$.
It should be noted that throughout this paper we will refer to constructing such a function as `interpolation', in keeping with the literature, although $\mathcal{F}$ may not be an exact fit to the data. 

\cite{Hoppe1992} approached this problem by imagining $\mathcal{F}(\vec{x})$ as the orthogonal projection of $\vec{x}$ onto the tangent plane of the surface at the nearest neighbour in the point cloud, $\vec{x}_j$:
\[ \mathcal{F}(\vec{x}) = (\vec{x}-\vec{x}_j)\cdot\vec{n}_j, \]
where each unit normal $\vec{n}_j$ defining the tangent plane at $\vec{x}_j$ is estimated by principal component analysis. 

Many other algorithms take advantage of surface normal information. 
\cite{Carr2001} constructed the scalar field $\mathcal{F}$ as a regularised radial basis function interpolant, fitted to an augmented point cloud with off-surface points added along normal directions to avoid the trivial (zero-everywhere) interpolant. 
Similarly, \cite{Ohtake2003} interpolated augmented point clouds with many low degree polynomials in a partition of unity. 
The partition of unity `blended' together local interpolants fit to small subsets of the data. 

Another implicit method proposed was Poisson surface reconstruction \citep{Kazhdan2006}, which sought not to enforce values at points ($\mathcal{F}(\vec{x_j}) = 0$), but rather to enforce normals: $\nabla \mathcal{F}(\vec{x}_j) = \vec{n}_j$.
This formulation was solved as a Poisson problem, and produced smooth surfaces robust to noisy data \citep{Kazhdan2006}.
Fitting only to the normals was found by \cite{Kazhdan2013} to produce over-smoothed surfaces. 
A variant called screened Poisson surface reconstruction was proposed \citep{Kazhdan2013} that balances the fit to the normals with the fit to the points.
Poisson reconstruction has recently been successful in reconstructing soy bean plants from point clouds for analysis of phenotype development \citep{Zhu2020}.

We choose polyharmonic spline radial basis function (RBF) interpolation to construct the scalar field $\mathcal{F}$ for its flexibility and smoothing properties as demonstrated by \cite{Carr2001}. 

\subsection{Noisy data}
Figure \ref{fig:noisyPC} shows a noisy dataset typical of those generated by the Artec scanner we use. 
Clearly, we require interpolation methods that are robust to (hopefully small) inaccuracies in the data. 
In least squares approximation problems, a common extension of the cost function involves including a regularisation term.
For our purposes, this extra term can control the smoothness of the interpolant.
This formulation is sometimes referred to as ridge regression, and was discussed at length by \cite{wahba1990spline}.

One use of ridge regression, which has been extensively studied \citep{wahba1990spline} and applied to noisy range data by \cite{Carr2004}, is spline smoothing. 
The success of this method hinges on the use of polyharmonic splines, which are radial basis functions minimising a penalty functional analogous to `bending energy'. 
Smoothing spline interpolants are found by minimising a cost function that depends on both the error of the interpolation, and the bending energy penalty.
When the error of the interpolant dominates the cost function, the surface may be rough and overfit the data. 
Conversely, when the bending energy penalty dominates, the reconstructed surfaces tend to more primitive forms, and only loosely resemble the data.
\cite{wahba1990spline} discussed ways to choose this parameter, for example by generalised cross-validation. 
We choose these splines for ease of implementation and extensibility to higher orders of continuity.

\subsection{Reducing complexity}
The usual approach to an RBF interpolation problem with $N$ points involves the solution of a large linear system of dimension $N\times N$, which is known to require $\bigO(N^3)$ floating point operations to solve using a typical direct method.
Such an approach, then, is clearly infeasible for a large number of points (\textcolor{revisedText}{e.g.,} $\bigO(10^5)$). 
Several methods have been proposed to reduce the computational complexity of solving the interpolation problem, many of which involve introducing locality to the interpolant.

\cite{wendland1995piecewise} introduced a class of RBFs with compact support to introduce sparsity (and therefore enable the use of fast sparse solvers), although this approach requires tuning the radius of support. 
\cite{Beatson2000} detailed the use of domain decomposition methods, in which solutions to smaller, local interpolation problems are used to precondition the global interpolation problem. 
When coupled with a fast method for computing the action of an interpolation matrix on a vector, this method enables efficient solutions to large problems; however,
fast evaluators (such as those by Beatson et al. \citep{Beatson1992,Beatson1997}) are technical to implement and RBF-specific.

The localisation method we consider is the partition of unity method (PUM) \textcolor{revisedText}{\citep{Franke1977,Franke1982,Babuska1997}. }
The PUM works by dividing the domain of the data into smaller, more manageable subdomains. 
The subdomains need not be distinct, but their union should cover the domain, and they should each be small enough to contain no more datapoints than is feasible for direct solution methods.
Figure \ref{fig:PUM} shows an example of a point cloud, partitioned into spherical volumes.
We fit local interpolants to the data in the subdomains, independent of one another.
It has been \textcolor{revisedText}{shown} that the PUM reduces the complexity of the interpolation problem on $N$ points to $\bigO(N)$ operations, under some mild assumptions on the data structures used in the implementation \citep{Wendland2002,Wendland2006}.
\textcolor{revisedText}{This RBF-PUM has been successfully applied to many problems, for example in the meshless discretisation of partial differential equations \citep{Cavoretto2019a,Heryudono2016}, and the characterisation of dynamical systems \citep{Cavoretto2016}.
A thorough discussion of the RBF-PUM in the context of multivariate interpolation can be found in \citep{Cavoretto2015}.}
Applications of the PUM with other interpolation methods can be found in the literature, such as the work by \cite{Ohtake2003}, using implicit polynomials as local interpolants. 

In this work, we propose a new technique based on implicit \textcolor{revisedText}{RBF-PUM interpolation} for reconstructing thin leaf surfaces, suitable for large, noisy datasets.
Our methods are most similar to the approach by \cite{Tobor2004} and \cite{Xiaojun2005}, who employ both RBF interpolants and the PUM to implicitly reconstruct closed surfaces.
We will show that the RBF-PUM, with polyharmonic smoothing splines, can reliably reconstruct smooth, thin surfaces from scanned point clouds.

The remainder of this paper is organised as follows: in section \ref{s:pre}, we detail the steps to prepare the dataset for fitting the implicit interpolant. 
Next, in section \ref{s:fit}, we discuss fitting the RBF-PUM interpolant.
In section \ref{s:post}, we focus on extracting a mesh representation of the implicitly defined surface reconstruction.
We present the results and discussion of our reconstructions of capsicum plant surfaces in section \ref{s:results}, and in section \ref{s:conclusion}, the outcome of our work is summarised.
\textcolor{revisedText}{
We have made our MATLAB implementation of the method available at: \url{https://github.com/rwhebell/ThinLeafReconstruction}.
}
\section{Preprocessing} \label{s:pre}

\subsection{Cleaning the data}
\begin{figure}[h]
	\centering
	\begin{subfigure}{0.4\linewidth}
	\includegraphics[width=\linewidth]{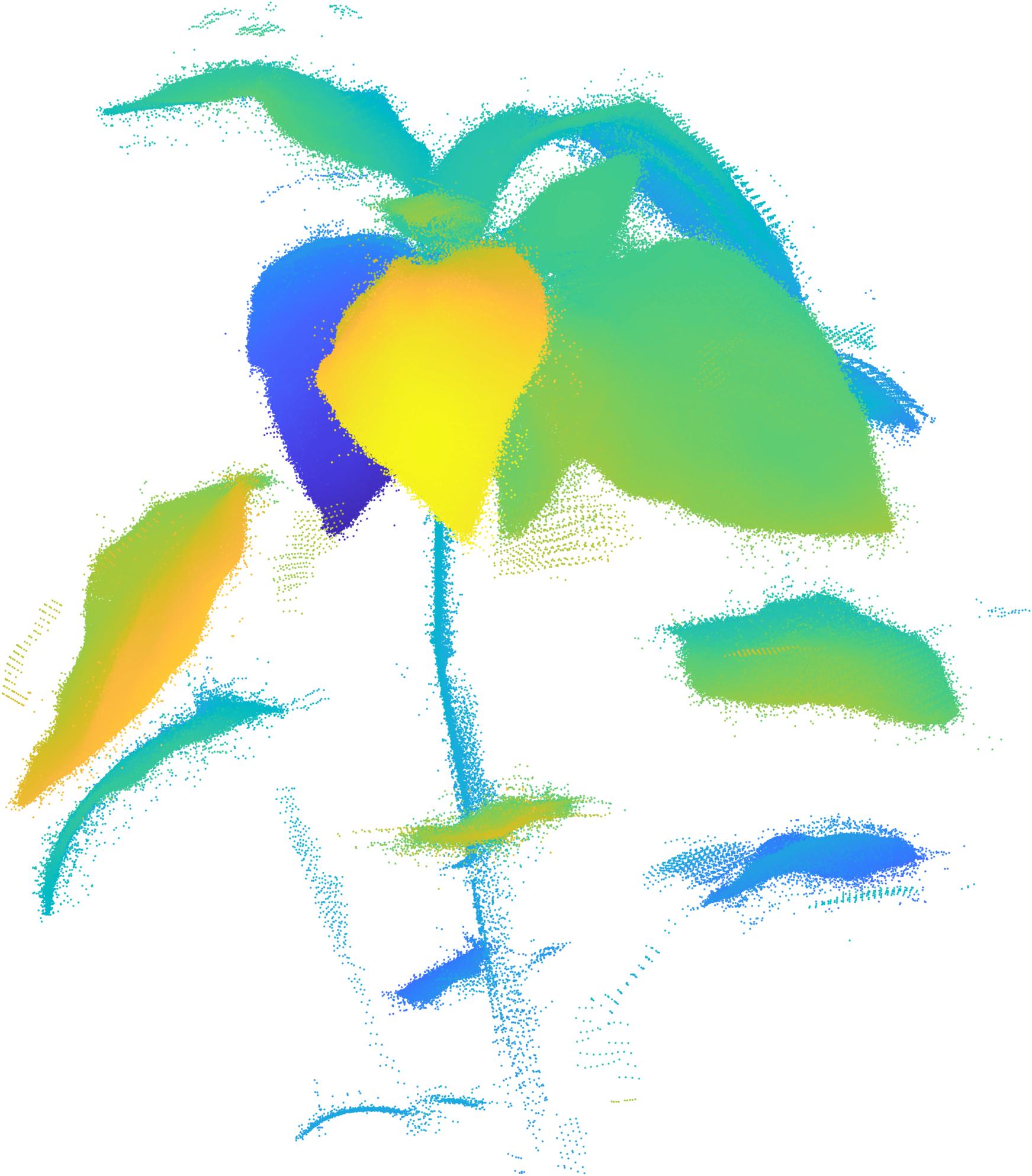}
	\end{subfigure}
	\hspace{12pt}
	\begin{subfigure}{0.4\linewidth}
	\includegraphics[width=\linewidth]{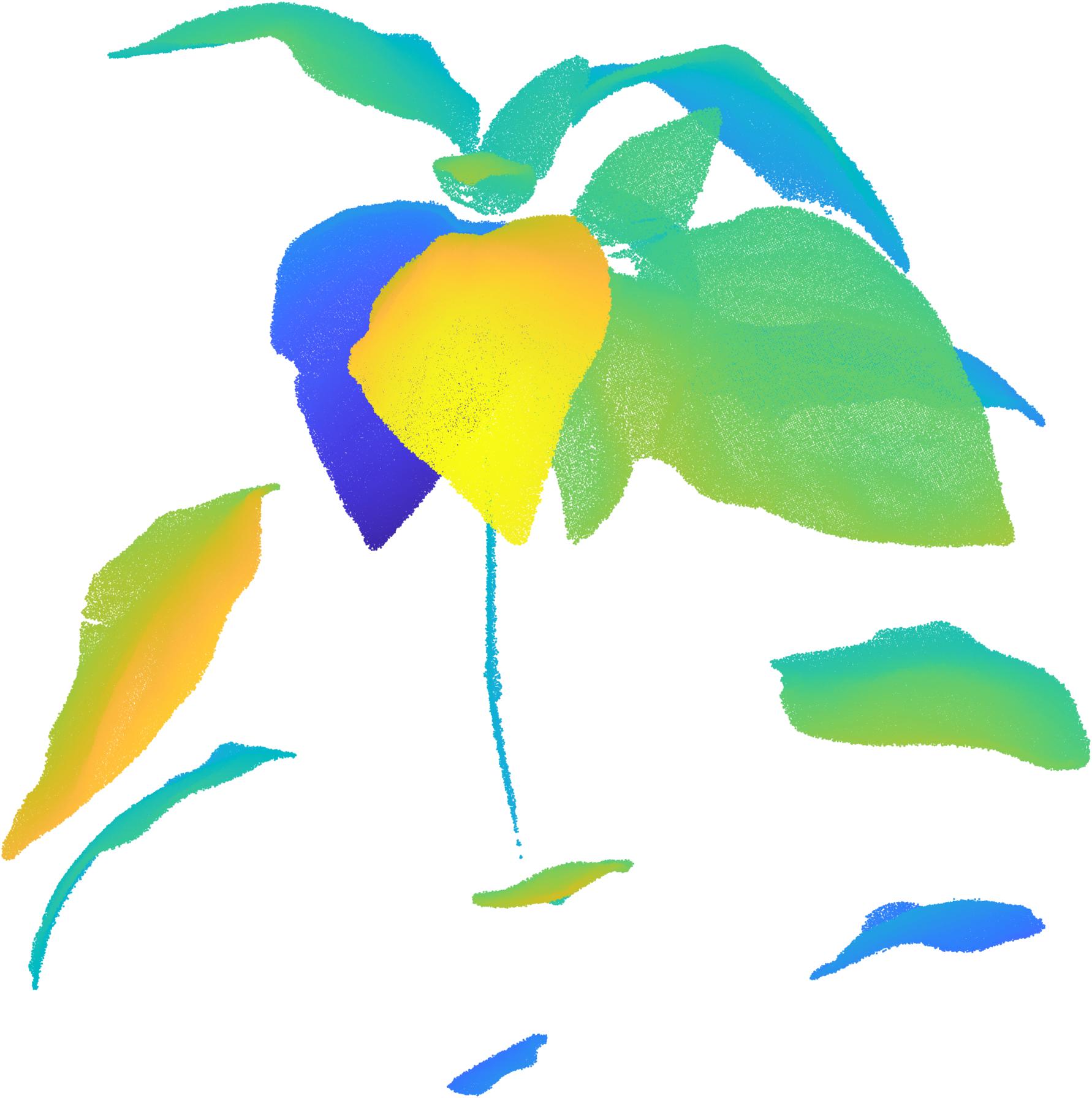}
	\end{subfigure}
	\caption{Scanned point cloud of a capsicum plant before (left) and after (right) preprocessing. Coloured by $x$-component for clarity.}
	\label{fig:noisyPC}
\end{figure}

Figure \ref{fig:noisyPC} shows a raw point cloud, scanned from a capsicum plant. 
Evidently, some outliers are present in the dataset, possibly due to plant movement during scanning, or poor registration (the matching of scans from different viewpoints).
Before we ready the dataset for interpolation, we clean it by removing outliers and downsampling to ensure a more regular point distribution.

The outlier detection we use is based on the average distance to neighbours of a point. 
We typically consider the average (Euclidean) distance to the 50 nearest neighbours of each point in the point cloud. 
Any point with an average greater than a specified threshold is labelled an outlier and excluded from the dataset.
For example, we use a threshold of 0.15 standard deviations for the capsicum plant shown in Figure \ref{fig:noisyPC}.
Once the point cloud is denoised, we downsample it using a grid average method \citep{MATLAB2020a} for a few reasons.
Firstly, we find the sampling density of the Artec Eva scanner more than sufficient for our purposes, so downsampling can offer computational speedup without adversely affecting the quality of the reconstruction.
Secondly, the condition number of the matrix arising from the RBF-PUM interpolation problem increases exponentially with the minimum point separation distance \citep{fasshauer2007meshfree} -- grid average downsampling combines very close-together points, increasing the minimum separation distance.
Lastly, the averaging may partially mitigate the effects of noise in the dataset.

\subsection{Estimating surface normals}
A closed, orientable surface $\sigma$ can be thought of as the zero level set of its signed-distance function. 
That is, if we have a scalar field $D(\vec{x})$ that is negative on one side of the surface, positive on the other, and satisfies:
\begin{align*}
|D(\vec{x})| = \min_{\vec{y}\in\sigma} \norm{\vec{y} - \vec{x}}{2},
\end{align*}
then we can recover the surface as its zero level set:
\begin{align*}
\sigma = \{ \vec{x} : D(\vec{x}) = 0 \}.
\end{align*}
We can fit interpolants that mimic signed distance functions and recover (orientable) surfaces as their zero level set, regardless of the surface's topology.

To fit such an interpolant given only points on the surface, we must add artificial points.
With only on-surface points, we would be fitting an interpolant to only zero-valued data, which would lead to the trivial solution that the interpolant is simply zero everywhere \citep{Carr2001}. 
Adding off-surface points with non-zero values to the dataset prevents this trivial solution.
It is typical to approximate surface normals from the point cloud and add off-surface points along these normals \citep{Hoppe1992}.
We use \citeauthor{MATLAB2020a}'s in-built method \texttt{pcnormals} to approximate surface normals at the datapoints. 
The algorithm is based on principal component analysis of a point's $k$ nearest neighbours (we typically use $k=50$ for our applications), and is described by \cite{Hoppe1992}.
We add the off-surface points $\vec{x}_j \pm L\vec{n}_j$. The scalar $L$ is on the order of the grid spacing used for downsampling.

Now that we have surface normals, we are faced with the problem of their orientation.
Clearly, $-\vec{n}_j$ would also be a valid choice for the normal if we were not concerned with the orientation of the surface. 
We must align the approximated normals by flipping their direction to define a consistent orientation, however it is not clear how we should traverse the surface to do so.
A solution to this problem is given by \cite{Hoppe1992}.
First, consider the weighted graph $G$, with vertices corresponding to the points $\vec{x}_j$ in the point cloud. 
If two points, $\vec{x}_i$ and $\vec{x}_j$ are sufficiently close, connect their nodes with an edge with weight $G_{ij}$, defined by:
\begin{equation*}
    G_{ij} = 1 - | \vec{n}_i \cdot \vec{n}_j |.
\end{equation*}
This way, a low edge weight $G_{ij}$ indicates a small change in the surface normal between points $\vec{x}_i$ and $\vec{x}_j$.
We consider two points to be sufficiently close if the Euclidean distance between them is less than some threshold, $\epsilon$.
Note that this condition means $G$ may not necessarily be connected.

\begin{figure}[h]
    \centering
    \begin{subfigure}[c]{0.45\linewidth}
        \includegraphics[width=\linewidth]{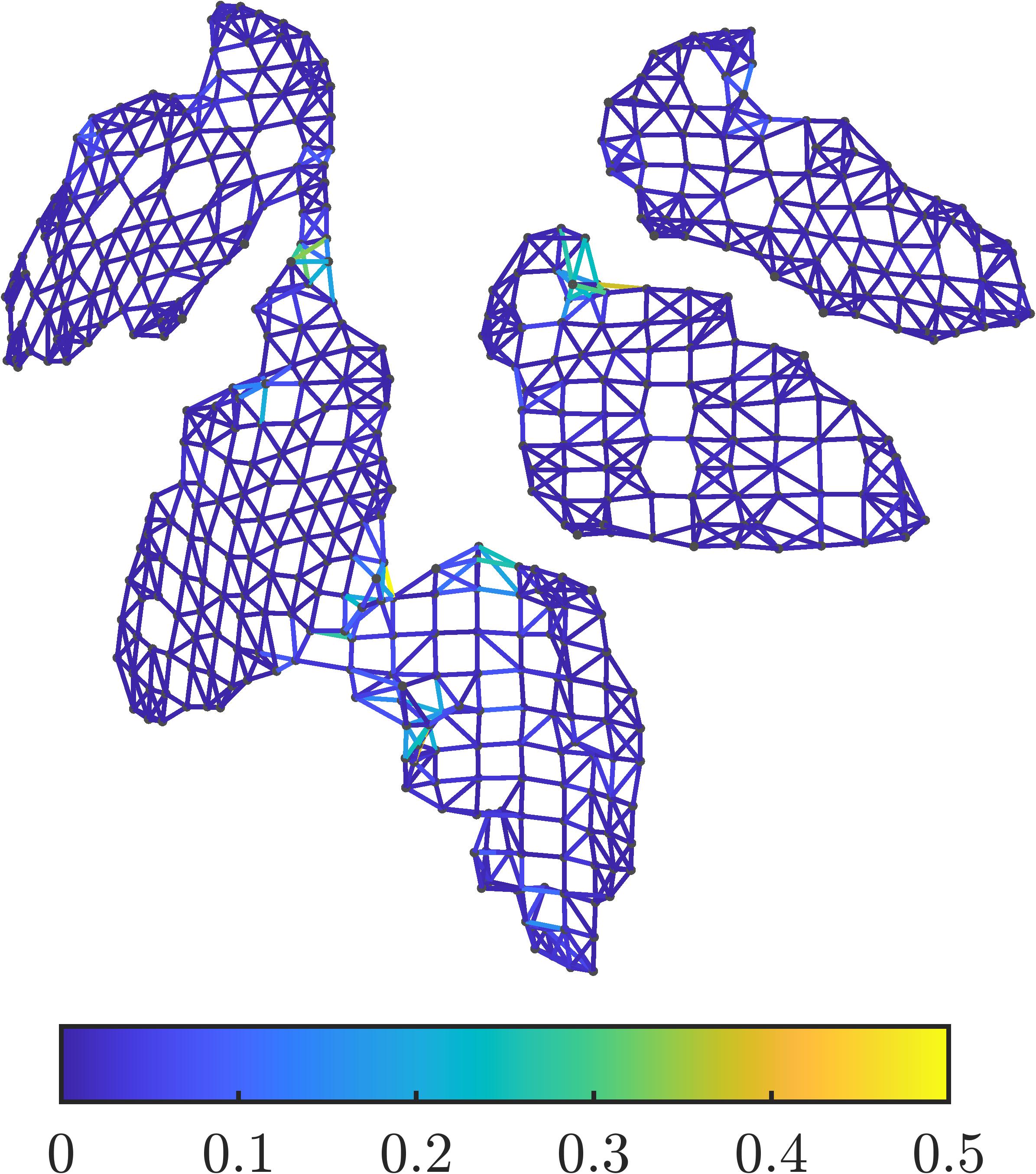}
    \end{subfigure}
    \hspace{12pt}
    \begin{subfigure}[c]{0.45\linewidth}
        \includegraphics[width=\linewidth]{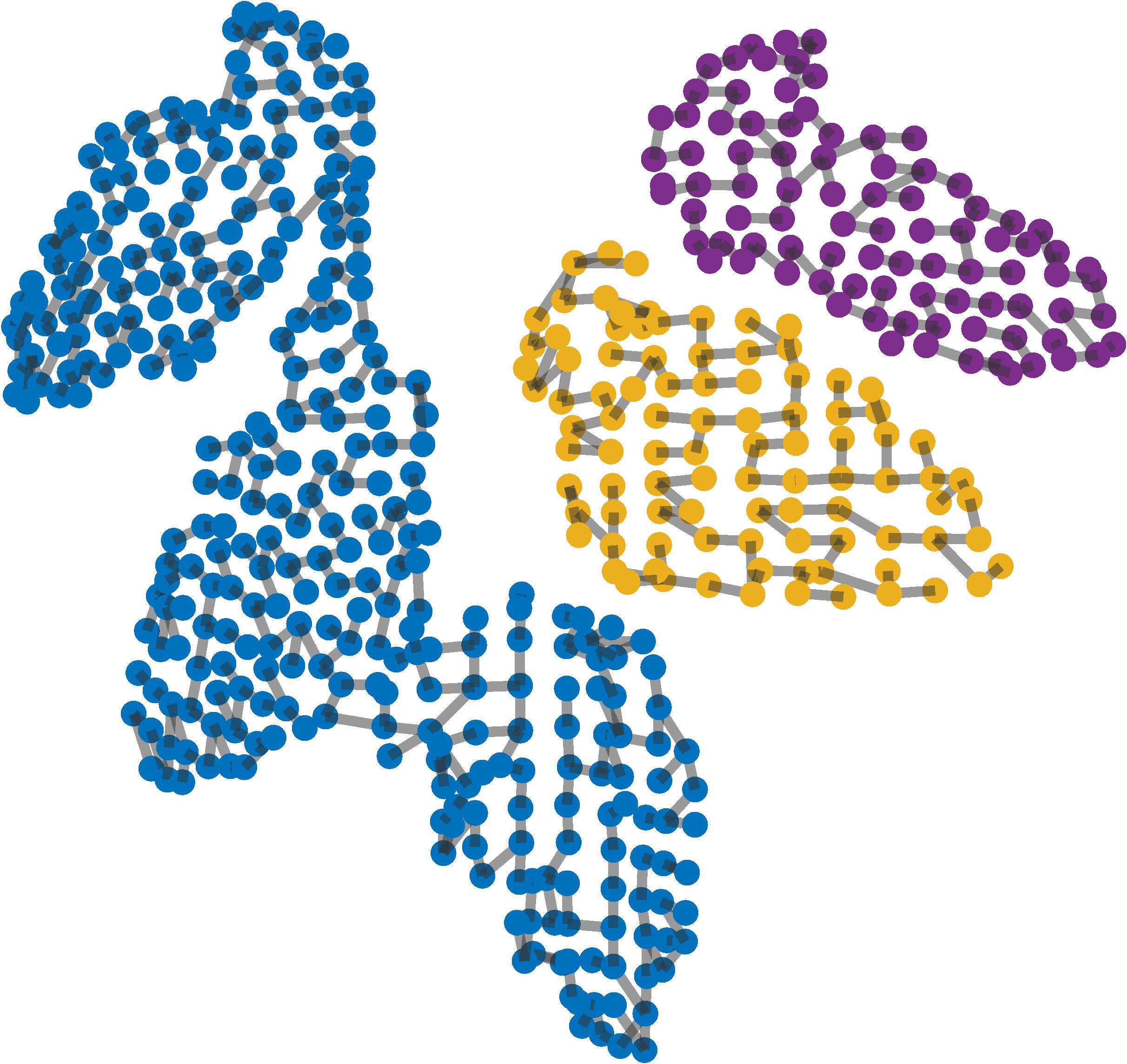} \vspace{0.8cm}
    \end{subfigure}
    \caption{Orienting normals for tomato leaves. Left: the weighted graph $G$, with edges coloured by weight. A high weight indicates a large change in normal. Right: a minimal spanning forest for the graph $G$, with nodes coloured to indicate distinct components.}
    \label{fig:tomatograph}
\end{figure}

We could traverse $G$ in search of neighbouring points with inconsistent normals (varying by more than $90\degree$), however we are likely to encounter cycles in the graph, or contradictory orientations. 
Instead, we find an approximate minimal spanning tree for each connected component of $G$ (a minimal spanning forest).
We then traverse these trees systematically, flipping normals where necessary.
We follow \cite{Hoppe1992} and use a breadth-first search.
This way we can avoid comparing points whose normals are dissimilar, where it may not be clear which point's normal orientation is `correct'. 
Figure \ref{fig:tomatograph} shows a graph and minimal spanning forest formed in this way for a downsampled point cloud of some tomato leaves.
Note that for large point clouds, it is sufficient to orient normals of a downsampled copy, and use this to orient the normals of the full point cloud. 
\textcolor{revisedText}{We have given pseudocode for this method in \ref{alg:orientNormals}.}

A useful side effect of this algorithm is that it clusters the points -- the components of the minimal spanning forest are usually distinct leaves, or at least groups of leaves close together.
This information could be used to reconstruct the leaves separately (in parallel).
Furthermore, in future droplet spray models, we may wish to extend existing capabilities \citep{mayo2015simulating,Dorr2014,Dorr2016} by varying certain wettability properties not just between plants, but also between leaves, making this clustering information very valuable.
The idea could also be extended to infer some chemical properties of the leaf by recognising its growth stage from its physical characteristics.

With the point cloud now equipped with oriented normals, we can now form implicit RBF-PUM interpolants \textcolor{revisedText}{with the set of surface points $\{ \vec{x}_1, \dots \vec{x}_N \}$, their values, $f_j = 0$, and the off-surface constraints:
\begin{gather*}
    \vec{x}_{N+j} = \vec{x}_j + L\vec{n}_j, \quad
    f_{N+j} = L,
    \\
    \vec{x}_{2N+j} = \vec{x}_j - L\vec{n}_j, \quad
    f_{2N+j} = -L.
\end{gather*}
}

\section{Interpolation} \label{s:fit}
Radial basis function (RBF) methods are powerful tools for interpolation which make very few assumptions about the structure of the data. 
We build interpolants as weighted sums of rotationally invariant basis functions, shifted to datasites, augmented with polynomial terms, in the form:
\begin{align*}
    \mathcal{F}(\vec{x}) = \sum_{j=1}^N \lambda_j \phi( \| \vec{x} - \vec{x}_j \|_2 )
    + \sum_{k=1}^{n} a_k p_k(\vec{x}),
    \quad
    \vec{x} \in \Omega,
\end{align*}
where $\vec{x}_j$ are the centres, $\phi(r)$ is the radial basis function, $p_k$ are the polynomials, $\lambda_i, a_k$ are weights to be determined, and $n$ is the dimension of the space of polynomials of degree at most $m-1$.

\subsection{Polyharmonic spline radial basis functions}\label{s:splines}
Polyharmonic splines are a class of radial functions on $\mathbb{R}^d$ that minimise a particular measure of the curvature of an interpolant. 
We will follow the derivation in the monograph by \cite{wahba1990spline}, where more detail and background can be found.

The `measure of curvature' we consider is the thin-plate penalty functional of order $m$ on $\mathbb{R}^d$: 
\begin{equation} \label{energy}
	J_m^d(\mathcal{F}) = 
	\sum_{\alpha_1+\dots+\alpha_s=m} \:
	\frac{m!}{\alpha_1! \dots \alpha_d!}
	\int_{\mathbb{R}^d}
	\left(
	\frac{ \partial^m \mathcal{F} }{ \partial x_1^{\alpha_1}\dots\partial x_d^{\alpha_d} }
	\right)^2
	\: \mathrm{d}\vec{x},
\end{equation}
which is the sum of squared $m$-th order mixed partial derivatives of $\mathcal{F}$.
For example, the second order ($m=2$) thin-plate penalty functional in $\mathbb{R}^2$ is:
\begin{equation*} 
	J_2^2(\mathcal{F}) = 
	\int_{\mathbb{R}^2}
	\left(
	\frac{ \partial^2 \mathcal{F} }{ \partial x_1^2 }
	\right)^2
	+
	2 \left(
	\frac{ \partial^2 \mathcal{F} }{ \partial x_1 x_2 }
	\right)^2
	+
	\left(
	\frac{ \partial^2 \mathcal{F} }{ \partial x_2^2 }
	\right)^2
	\: \mathrm{d}x_1 \: \mathrm{d}x_2.
\end{equation*}
We restrict ourselves to interpolants in the space containing distributions on $\mathbb{R}^d$ with square-integrable $m$-th derivatives, so that the penalty functional \pref{energy} is defined.
Furthermore, for reasons discussed by \cite{wahba1990spline}, we must have $2m-d>0$.

\cite{duchon1977splines} showed that the basis functions of an interpolant minimising this penalty functional are Green's functions for the $m$-iterated Laplacian, $\Delta^m$ \citep{wahba1990spline}, namely:
\begin{align*}
	E_m(\vec{x},\vec{x}_j) &\coloneqq E(\norm{\vec{x}-\vec{x}_j}{2}), \text{ where}
	\\
	E(\tau) &= 
	\begin{cases}
		\theta_m^d |\tau|^{2m-d} \log|\tau|, & \text{$2m-d$ even}, \\
		\theta_m^d |\tau|^{2m-d}, & \text{otherwise},
	\end{cases}
	\\
	\text{with} \quad \theta_m^d &= 
	\begin{cases}
	\frac{(-1)^{d/2+1+m}}{ 2^{2m-1} \pi^{d/2} (m-1)! (m-d/2)! },
	& \text{$2m-d$ even}, \\
	\frac{\Gamma(d/2-m)}{ 2^{2m} \pi^{d/2} (m-1)! }, 
	& \text{otherwise}.
	\end{cases}
\end{align*}
Formally, the functions $E_m(\cdot,\vec{x}_j)$ solve the polyharmonic equation:
\begin{equation*}
\Delta^m E_m( \cdot, \vec{x}_j ) = \delta(\vec{x}_j),
\end{equation*}
where $\delta(\vec{x}_j)$ is the Dirac delta function, so that
\begin{equation*}
\Delta^m \mathcal{F}(\vec{x}) = 0, \quad
\text{ for } \vec{x}\not=\vec{x}_j, \quad 
j=1,\dots,N.
\end{equation*}
Thus, an interpolant minimising the penalty functional \pref{energy} has the form:
\begin{equation*}
	\mathcal{F}(\vec{x}) = 
	\sum_{j=1}^N \lambda_j E_m(\vec{x},\vec{x}_j) +
	\sum_{k=1}^{n} a_k p_k(\vec{x}),
\end{equation*}
where the polynomials $p_k$ span the $n = \binom{d+m-1}{d}$-dimensional space of polynomials of total degree at most $m-1$. 

The weights $\lambda_j$ and $a_k$ are determined by enforcing exact interpolation conditions:
\begin{equation*}
    \mathcal{F}(\vec{x}_j) = f_j,
\end{equation*}
imposing the additional orthogonality conditions \citep{wahba1990spline}:
\begin{align*}
    \sum_{j=1}^N \lambda_j p_k( \vec{x}_j ) = 0,
\end{align*}
and solving the resulting linear system:
\begin{equation} \label{polyInterp}
    \begin{gathered}
    \begin{bmatrix}
    	\vec{A} & \vec{P} \\
    	\vec{P}^T & \vec{0}_{n\times n}
    \end{bmatrix}
    \begin{bmatrix}
    	\vec{\lambda} \\ \vec{a}
    \end{bmatrix}
    =
    \begin{bmatrix}
    	\vec{f} \\ \vec{0}_{n\times 1}
    \end{bmatrix},
    \\
    A_{ij} = \phi( \norm{\vec{x}_i - \vec{x}_j}{2} ), \quad
    i,j = 1,\dots,N,
    \\
    P_{ik} = p_k( \vec{x}_i ), \quad
    i = 1,\dots,N, \: k = 1,\dots,n,
    \\
    \vec{f}=[f_1,\dots,f_N]^T. 
    \end{gathered}
\end{equation}
It can be shown that the polyharmonic spline of order $m$ on $\mathbb{R}^d$ is strictly conditionally positive definite of order $m$ \citep{duchon1977splines,wahba1990spline,fasshauer2007meshfree}.
That is, the interpolation conditions together with the orthogonality conditions have a unique solution (if the $\vec{x}_j$ are such that least squares regression on $p_k$ is unique).

Now that we have basis functions that form exact interpolants with minimal curvature, a natural extension is fitting approximate interpolants by adding a `curvature penalty'.
To fit approximate interpolants, we relax the exact interpolation conditions, and instead minimise the objective function:
\begin{equation} \label{costFn}
\frac{1}{N} \sum_{j=1}^N ( u_j - \mathcal{F}(\vec{x}_j) )^2
+ \rho J_m^d(\mathcal{F}),
\end{equation}
where $J_m^d$ is the penalty functional defined in equation \pref{energy} and $\rho$ is the smoothing parameter.
As mentioned in the introduction, this is known in the statistics literature as ridge regression \citep{wahba1990spline}.

The parameters $\vec{\lambda} = [\lambda_1,\dots,\lambda_N]^T$ and $\vec{a}=[a_1,\dots,a_n]^T$ can then be shown \citep{wahba1990spline} to solve the modified linear system:
\begin{gather} \label{smoothAlambdab}
    \begin{gathered}
	\begin{bmatrix}
		\vec{A} + (\rho N / \theta_m^d) \vec{I} & \vec{P} \\
		\vec{P}^T & \vec{0}_{n\times n}
	\end{bmatrix}
	\begin{bmatrix}
	\vec{\lambda} \\ \vec{a}
	\end{bmatrix}
	=
	\begin{bmatrix}
		\vec{f} \\ \vec{0}_{n\times 1}
	\end{bmatrix},
	\\
	A_{ij} = \phi( \norm{\vec{x}_i - \vec{x}_j}{2} ), \quad
    i,j = 1,\dots,N,
    \\
    P_{ik} = p_k( \vec{x}_i ), \quad
    i = 1,\dots,N, \: k = 1,\dots,n,
    \\
    \vec{f} = [f_1,\dots,f_N]^T,
    \end{gathered}
    \\[6pt]
	\phi( \tau ) = 
	\begin{cases}
		|\tau|^{2m-d} \log|\tau|, & \text{$2m-d$ even}, \\
		|\tau|^{2m-d}, & \text{otherwise},
	\end{cases}
\end{gather}
where $\vec{I}$ is the $N \times N$ identity matrix.

We now have a class of basis functions with smoothing properties for which we can guarantee the existence of unique solutions to our interpolation problem.
Table \ref{t:myfuncs} summarises some basis functions in the class.

\begin{table}[h]
	\centering
	\caption{Some polyharmonic spline radial basis functions.}
	\label{t:myfuncs}
	\begin{tabular}{cccc}
		\toprule
		Order ($m$) & Dimension ($d$) &        RBF ($\phi$)       & $1 / \theta_m^d$ \\ \midrule
		     2      &        2        & $r^2 \log r$ & $8\pi$ \\
		     2      &        3        &     $r$      & $-8\pi$ \\
		     3      &        3        &    $r^3$     & $96\pi$ \\ 
		     4      &        3        &    $r^5$     & $-2880\pi$ \\ \bottomrule
	\end{tabular}
\end{table}

\subsection{Generalised cross validation}\label{s:gcv}
\textcolor{revisedText}{
Generalised cross validation (GCV) is a method for estimating the optimal smoothing parameter $\rho$ \citep{golub1979generalized}.
GCV is frequently used to tune smoothing parameters, for example when fitting regularised cubic splines to point data \citep{bertolazzi2020point};  reconstructing terrain models from LiDAR data \citep{Chen2019}; and indeed in reconstructing wheat leaves \citep{Kempthorne2015}. 
The optimal smoothing parameter according to GCV minimises the objective function \citep{wahba1990spline}:
\begin{gather} \label{eq:gcv}
    V(\rho) = \frac{ N \norm{\left( \vec{I}-\vec{B}(\rho) \right) \vec{f} }{2}^2 }{ \text{trace}(\vec{I}-\vec{B}(\rho))^2 },
\end{gather}
where $\vec{f}$ is the vector of the sample values $f_j$, $\vec{I}$ is the $N\times N$ identity matrix, and $\vec{B}(\rho)$ is the influence matrix, which is defined such that:
\begin{gather*}
    \begin{bmatrix} \mathcal{F}(\vec{x}_1) \\ \vdots \\ \mathcal{F}(\vec{x}_N) \end{bmatrix}
    =
    \vec{B}(\rho) \vec{f}.
\end{gather*}
Letting the QR decomposition of the matrix $\vec{P}$ from equation \pref{smoothAlambdab} be:
\begin{gather*}
    \vec{P} = 
    \begin{bmatrix} \vec{Q}_1 & \vec{Q}_2 \end{bmatrix}
    \begin{bmatrix} \vec{R} \\ \vec{0}_{n\times1} \end{bmatrix},
\end{gather*}
we have the convenient representation for the matrix $\vec{I}-\vec{B}(\rho)$ due to \cite{wahba1990spline}:
\begin{gather*}
    \vec{I} - \vec{B}(\rho) = 
    \left( \rho N / \theta_m^d \right) \vec{Q}_2 \left( \vec{Q}_2^T \left( \vec{A} + (\rho N/\theta_m^d) \vec{I} \right) \vec{Q}_2 \right)^{-1} \vec{Q}_2^T.
\end{gather*}
}
We find an estimate for the smoothing parameter $\rho$ by minimising $V(\rho)$.
MATLAB's \citep{MATLAB2020a} \texttt{fminbnd} function is well suited to this, given a sensible region in which the smoothing parameter is expected to lie. 
For example, we search for $\rho$ between $10^{-6}$ and $10^{-1}$.
\subsection{Partition of unity method}
Direct interpolation of scattered data quickly becomes infeasible for large datasets. 
Notice that for a dataset with $N$ datasites, we require $\mathcal{O}(N^2)$ storage and $\mathcal{O}(N^3)$ operations for the solution to the linear system (\ref{polyInterp}). 
Several techniques exist to reduce the complexity of this interpolation problem.
We use the partition of unity method (PUM) for introducing locality to the interpolant, so that we may decouple the solution of the interpolation problem into many smaller problems.

The PUM was first discussed by Franke \citep{Franke1977}, and analysed by \cite{Wendland2002}. 
The core idea is to partition the domain into $M$ (possibly overlapping) subdomains. 
A local interpolant is then fit to each of these independently, and the global interpolant is a weighted sum of the local interpolants. 
The PUM is efficient, and does not sacrifice the exactness of the global interpolant, as we will discuss.

First, we partition the domain, $\Omega$, into $M$ overlapping sub-domains $\{\Omega_i\}_{i=1}^M$ such that the union of the subdomains contains the domain entirely. 
That is,
\begin{equation*}
\Omega \subseteq \bigcup_{i=1}^M \Omega_i.
\end{equation*}
We consider $M$ interpolants -- one for the data in each of the subdomains. 
Note that this method can lead to one datapoint being used to fit multiple interpolants. 
For efficiency, then, the overlap of neighbouring subdomains should be small. 
Figure \ref{fig:PUM} shows an example of a domain partitioned in this way. 
We use spherical (in three dimensions) subdomains, defined by a centre, $\vec{c}_i$, and a radius, $r_i$. 
The subdomain is then given by:
\begin{equation*}
\Omega_i := \left\{ \vec{x} \: : \: \frac{\norm{\vec{x}-\vec{c}_i}{2}}{r_i} \leq 1 \right\}.
\end{equation*}
We denote the interpolant of the data in $\Omega_i$ as $\mathcal{F}_i$.

Next, we require a way to formulate a global interpolant from these local interpolants. We do this by weighting the contribution of a patch by the distance to its centre. The weight functions we use are compactly supported Wendland functions \citep{wendland1995piecewise} such as:
\begin{equation*}
\varphi(r) =
	\begin{cases}
		(1-r)^4(4r+1), & \text{if } 0\leq r \leq 1 \\
		0, & \text{otherwise},
	\end{cases}
\end{equation*}
or:
\begin{align*}
\varphi(r) =
	\begin{cases}
		(1 - r)^6 (35r^2 + 18r + 3), & \text{if } 0\leq r \leq 1 \\
		0, & \text{otherwise},
	\end{cases}
\end{align*}
depending on the degree of continuity required. 

The weight functions must be normalised to sum to unity at any $\vec{x}\in\Omega$. First, we map the input to the unit sphere:
\begin{equation*}
\varphi_i(\vec{x}) =  \varphi \left( 
\frac{ \norm{\vec{x}-\vec{c}_i}{2} }{ r_i } \right),
\quad i = 1,\dots,M.
\end{equation*}
Next, we use Shepard's method \citep{shepard1968two} to normalise these functions:
\begin{equation*}
w_i(\vec{x}) = 
\frac{\varphi_i(\vec{x})}{ \sum_{k=1}^M \varphi_k(\vec{x}) }.
\end{equation*}
Finally, we construct the global interpolant:
\begin{align}
\mathcal{F}(\vec{x}) &= \sum_{i=1}^M w_i(\vec{x}) \mathcal{F}_i(\vec{x}) \nonumber 
\\
&= \sum_{i=1}^M w_i(\vec{x}) 
\left\{
\sum_{j=1}^{N_i} 
\lambda_j^i \phi \left( \norm{\vec{x}-\vec{x}_j^i}{2} \right) +
\sum_{k=1}^{n} a_k^i p_k( \vec{x} )
\right\}.
\end{align}

We note that if the local interpolants are exact, the global interpolant respects this exactness. To see this connection, consider a sample point $(\vec{x}_j,u_j)$. The local interpolants satisfy:
\begin{align*}
\mathcal{F}_i(\vec{x}_j^i) = f_j^i, \quad \vec{x}_j^i \in \Omega_i.
\end{align*}
Now, since the weight functions are normalised to sum to unity, and zero outside their respective subdomains, we have the global interpolant:
\begin{align*}
\mathcal{F}(\vec{x}_j^i) &= \sum_{i=1}^M w_i(\vec{x}_j^i) \mathcal{F}_i(\vec{x}_j^i)
= \sum_{i=1}^M w_i(\vec{x}_j^i) f_j^i
= f_j.
\end{align*}

Furthermore, the PUM is stable in the sense that the global error is bounded by the worst local error \citep{Wendland2006}. To see this, consider a function $g(\vec{x})$, and its PUM interpolant $\mathcal{F}(\vec{x})$, from samples $(\vec{x}_j,g(\vec{x}_j))$, and then take the absolute error to obtain:
\begin{align*}
\left| g(\vec{x}) - \mathcal{F}(\vec{x}) \right|
&=
\left| \sum_{i=1}^M w_i(\vec{x}) \left[ g(\vec{x}) - \mathcal{F}_i(\vec{x}) \right] \right|
\\
&\leq
\sum_{i=1}^M w_i(\vec{x}) \left| g(\vec{x}) - \mathcal{F}_i(\vec{x}) \right|
\\
&\leq 
\max_{1\leq i\leq M} \| g - \mathcal{F}_i \|_{ L_{\infty}(\Omega_i) }.
\end{align*}

\subsection{Choosing subdomains} \label{s:fit_choice}
If the data were somewhat evenly scattered, we could simply set the subdomains on a regular grid. 
For other types of data (point clouds, for example) we must be more careful about how we partition the domain.

We use an octree-like method similar to that described by \cite{Ohtake2003} to partition the domain in the case of unorganised point cloud data.
The method is given in psuedo-code in \textcolor{revisedText}{\ref{PUMoctree}}.
The core idea is to iteratively partition the domain into distinct cubes, then construct a spherical subdomain to cover each cube (we refer to these as \textit{covering spheres}). 
We introduce two new parameters for this method, $n_\text{min}$ and $n_\text{max}$. 
At each iteration, the cube whose covering sphere contains the most datasites is partitioned into eight sub-cubes, until no covering sphere contains more than $n_\text{max}$ datapoints.
Once these iterations are halted, we ensure that each subdomain has sufficiently many datapoints by expanding the radii of spheres containing fewer than $n_\text{min}$ datapoints.

\begin{figure}[h]
	\centering
	\begin{subfigure}[t]{0.31\linewidth}
		\includegraphics[width=\textwidth]{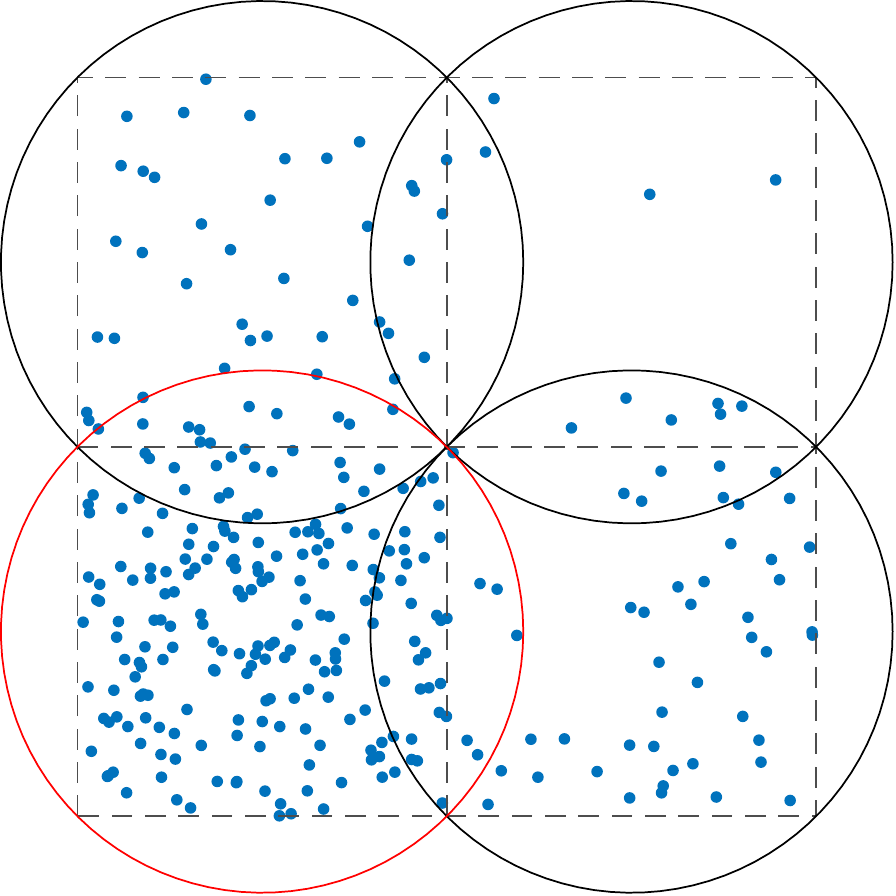}
		\caption{The bottom-left subdomain contains more than $n_\text{max}$ points, so it is replaced with four smaller subdomains.}
	\end{subfigure}
	~
	\begin{subfigure}[t]{0.31\linewidth}
		\includegraphics[width=\textwidth]{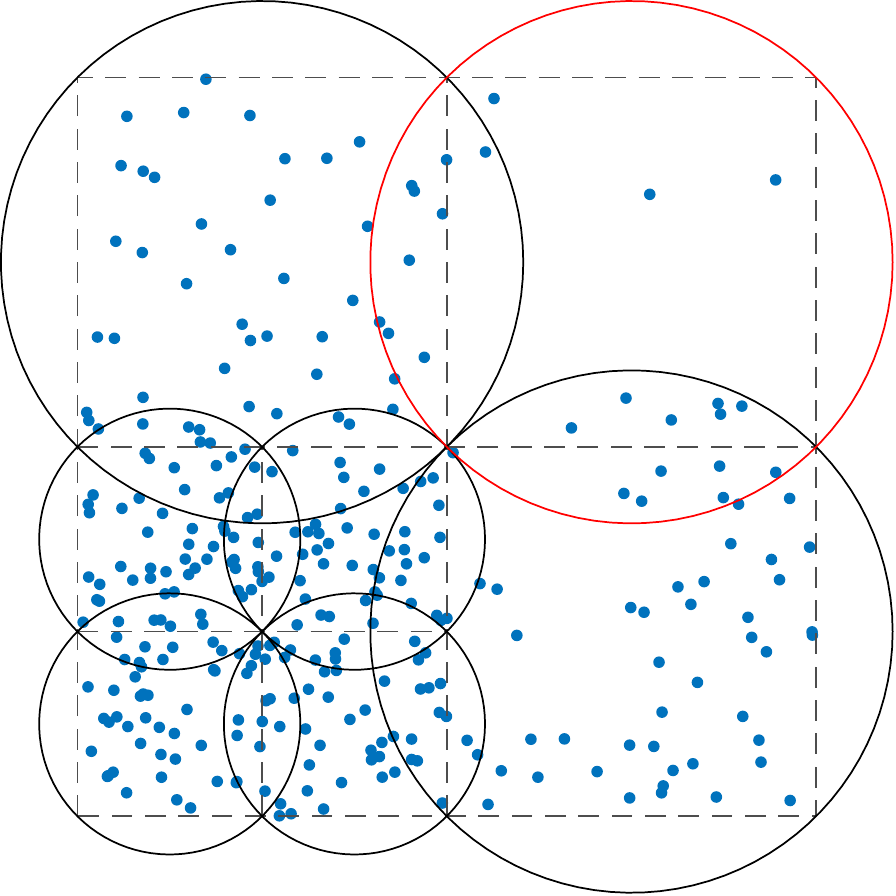}
		\caption{The top-right subdomain contains fewer than $n_\text{min}$ points, so it is expanded.}
	\end{subfigure}
	~
	\begin{subfigure}[t]{0.31\linewidth}
		\includegraphics[width=\textwidth]{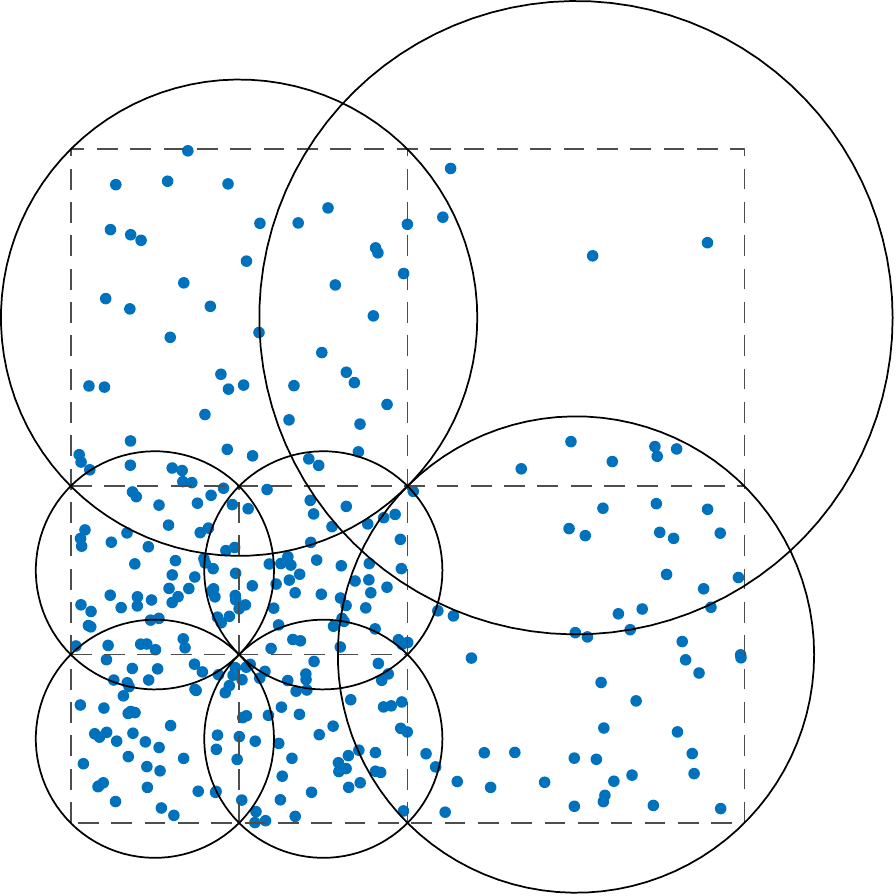}
		\caption{The final partition of the domain, satisfying $N_i \in [n_\text{min},n_\text{max}]$ for $i=1,\dots,M$.}
	\end{subfigure}
	\caption[Example of octree-like method for partitioning a domain.]{Example of an octree-like method for partitioning a domain. Sample points (blue) in a domain in $\mathbb{R}^2$, with varying point density. Black circles show the boundaries of the circular subdomains. Dashed grey lines show the underlying octree partitioning.}
	\label{fig:octree}
\end{figure}

\begin{figure}[h]
	\centering
	\begin{subfigure}{0.4\linewidth}
    	\includegraphics[width=\linewidth]{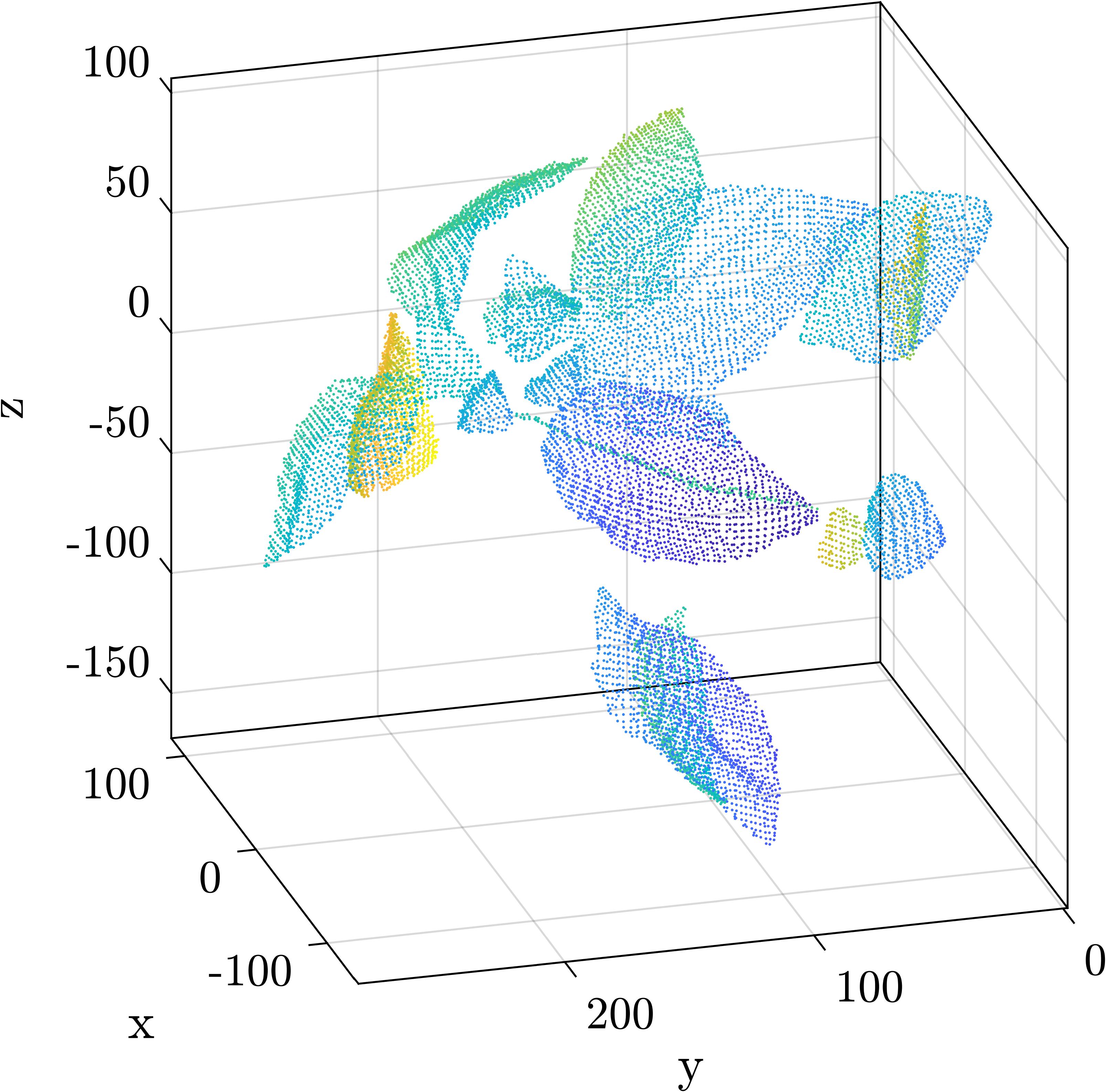}
	\end{subfigure}
	\hspace{12pt}
	\begin{subfigure}{0.4\linewidth}
    	\includegraphics[width=\linewidth]{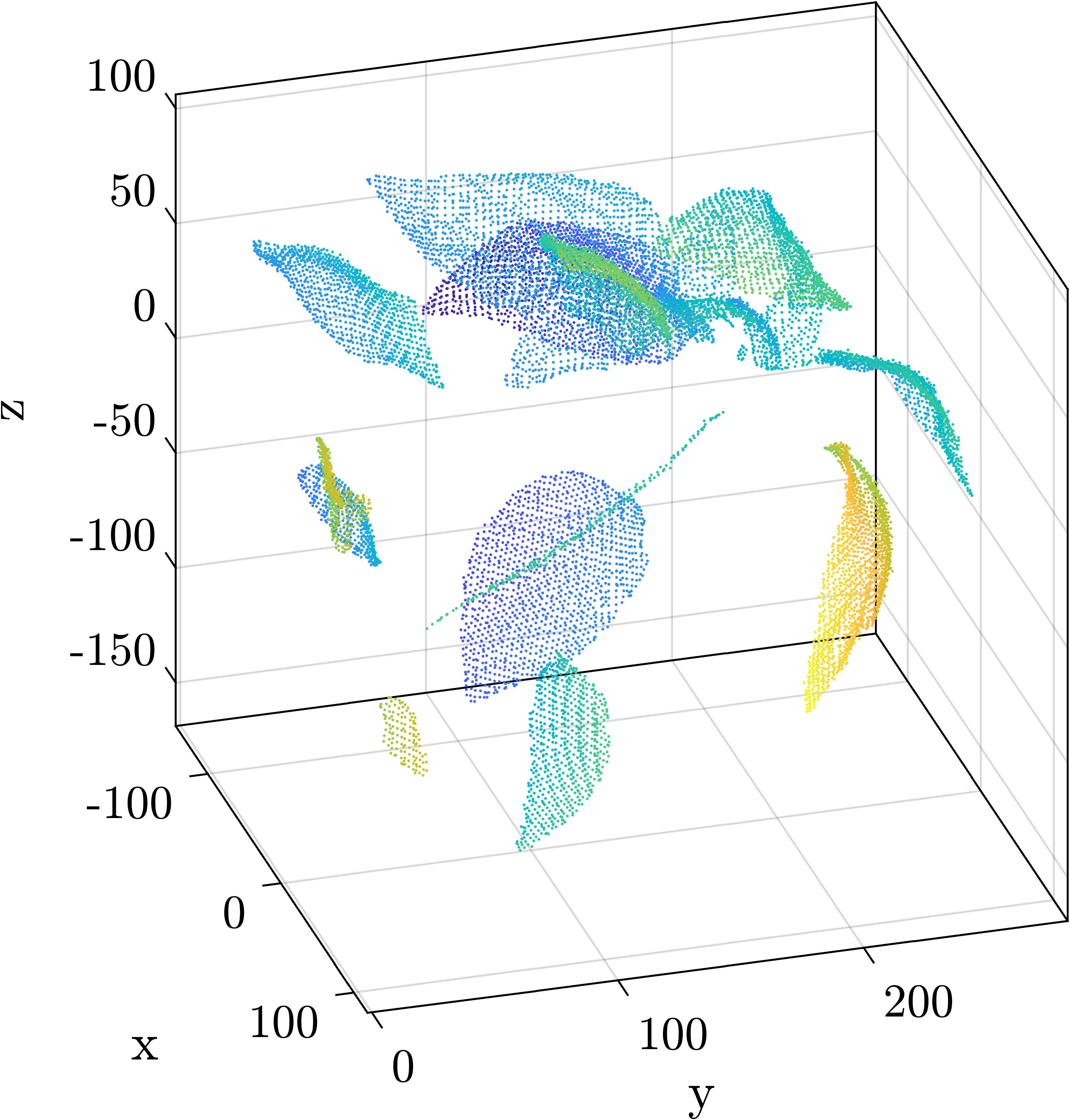}
	\end{subfigure}
	\\[12pt]
	\begin{subfigure}{0.4\linewidth}
	    \includegraphics[width=\linewidth]{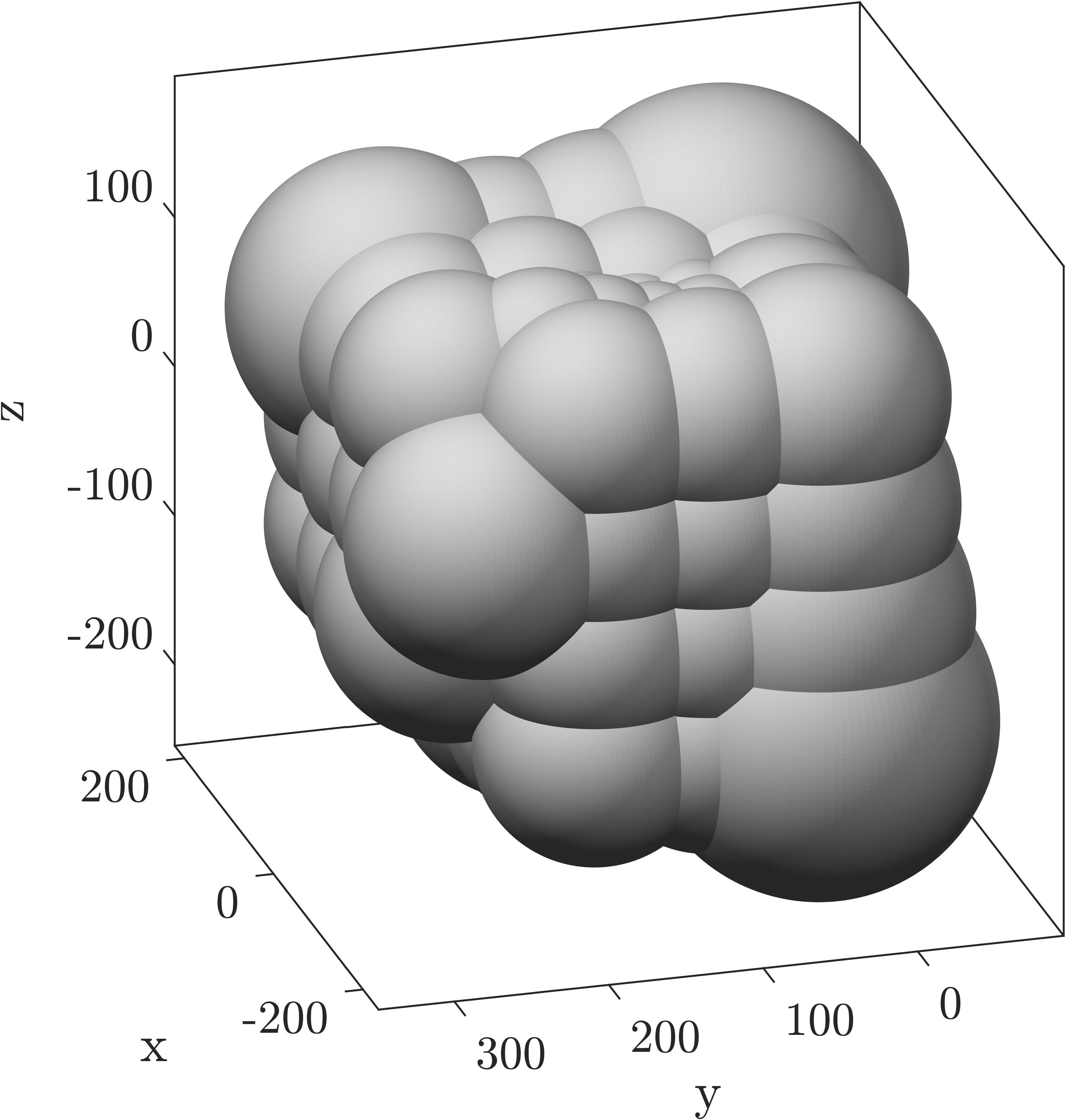}
	\end{subfigure}
	\hspace{12pt}
	\begin{subfigure}{0.4\linewidth}
	    \includegraphics[width=\linewidth]{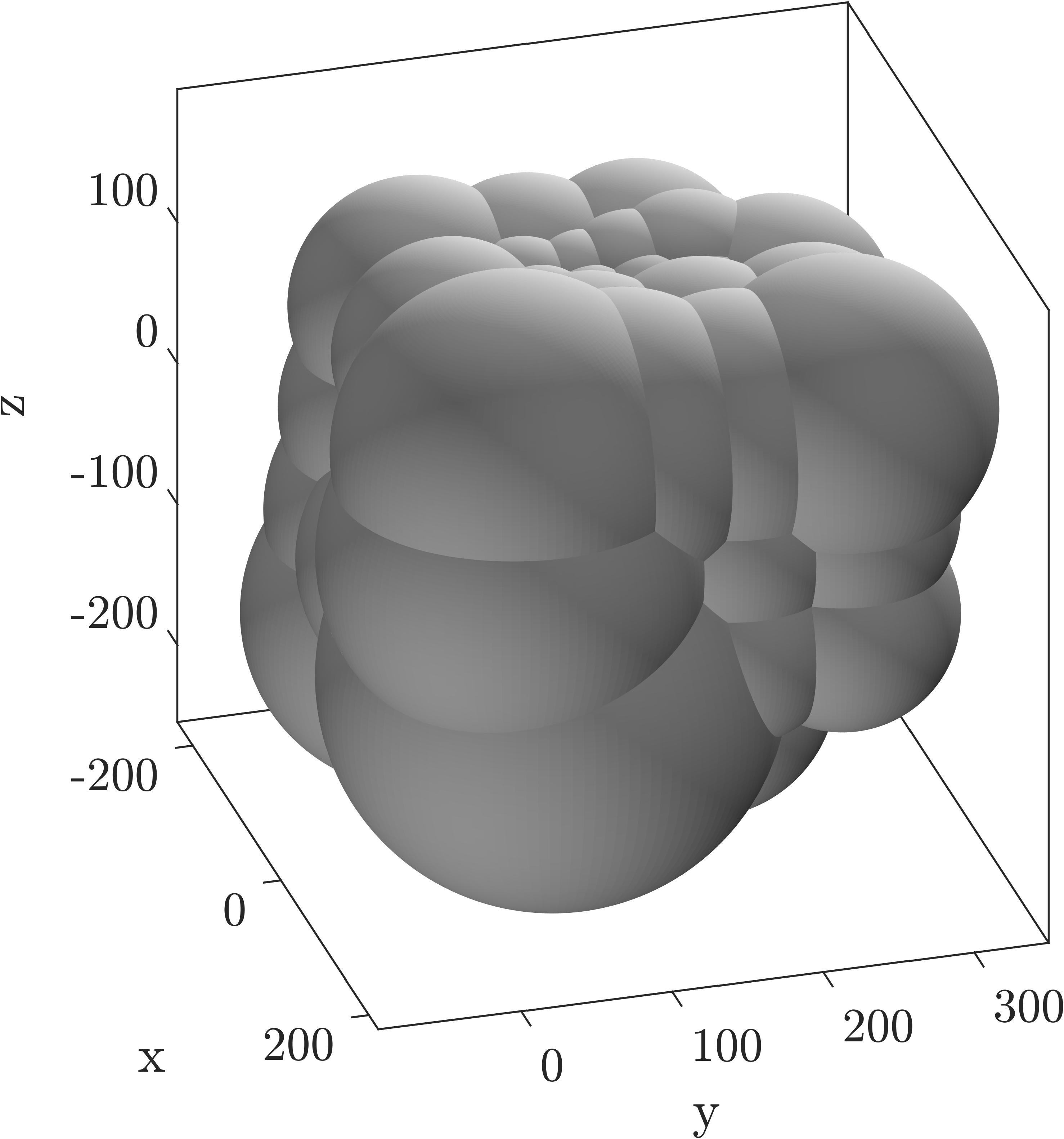}
	\end{subfigure}
	\caption{Top: point cloud of a capsicum plant from two viewpoints. Bottom: spherical subdomains for the partition of unity method. Octree-like domain partitioning ensures each subdomain has between $n_\text{min}$ and $n_\text{max}$ points.}
	\label{fig:PUM}
\end{figure}

Figure \ref{fig:octree} shows a simple example of the algorithm for 2D, variable density data.
The first step ensures that no subdomain contains too many points, and the second ensures that each subdomain contains enough points to produce a quality interpolant.
Figure \ref{fig:PUM} shows the result of an octree-like partitioning for a capsicum plant point cloud.

\textcolor{revisedText}{
To investigate the time complexity of building this octree-like partition, we recorded the runtime of the algorithm when applied to a point cloud with random subsets of points discarded.
We ran 10 trials for 20 different values of $N$, the remaining number of points.
A power law fit to a log-log plot of mean runtime against $N$ indicates a time complexity of approximately $\bigO\left(N^{1.3}\right)$, including the time taken to organise the points into a kd-tree data structure \citep{Bentley1975} for efficient spatial searching.
We note, however, that the time spent constructing the octree-like partition does not represent a significant portion ($<1\%$) of the overall time taken to fit the interpolant.
}

\subsection{Numerical solutions for RBF-PUM interpolants}\label{s:fit_num}
As usual, we enforce the interpolation conditions:
\begin{gather*}
	\mathcal{F}_i(\vec{x}_j^i) = f_j^i, 
	\quad \vec{x}_j^i \in \Omega_i, 
	\\ 
	i = 1,\dots, M, 
	\quad j = 1, \dots, N_i.
\end{gather*}
At first glance, this equation yields a large (albeit sparse) linear system for the weights. 
For example, an interpolation problem with two subdomains is solved by the linear system:
\begin{gather*}
\begin{bmatrix}
	\vec{A}(\Omega_1) & \vec{P}(\Omega_1) & & 
	\\
	\vec{P}(\Omega_1)^T & \vec{0} & &
	\\ 
	& & \vec{A}(\Omega_2) & \vec{P}(\Omega_2)
	\\
	& & \vec{P}(\Omega_2)^T & \vec{0}
\end{bmatrix}
\begin{bmatrix}
\vec{\lambda}^1 \\ \vec{a}^1 \\ \vec{\lambda}^2 \\ \vec{a}^2
\end{bmatrix}
=
\begin{bmatrix}
\vec{f}^1 \\ \vec{0} \\ \vec{f}^2 \\ \vec{0}
\end{bmatrix},
\\[5pt]
\vec{A}(\Omega_l)_{ij} = \phi\left(\norm{\vec{x}_i-\vec{x}_j}{2}\right), \quad
\vec{x}_i, \vec{x}_j \in \Omega_l,
\\[5pt]
\vec{P}(\Omega_l)_{ik} = p_k(\vec{x}_i), \quad
\vec{x}_i \in \Omega_l,
\\[5pt]
\vec{\lambda}^i = [ \lambda_1^i, \dots, \lambda_{N_i}^i ]^T,
\\[5pt]
\vec{a}^i = [ a_1^i, \dots, a_n^i ]^T,
\end{gather*}
where $\vec{f}^i$ is the vector of samples corresponding to the points $\vec{x}$ in the subdomain $\Omega_i$.

We decouple this system into $M$ smaller linear systems -- one for each subdomain -- since the local interpolants are independent. 
These smaller linear systems could theoretically be solved in parallel, providing even greater speed-up.
Furthermore, the total operation count for solving all of the local interpolation problems is $\bigO(M(N/M)^3)$. 
If we consider $N/M$ (the average number of datapoints per subdomain) to be constant, this reduces to $\bigO(N)$.
Note that the PUM runtime depends not only on the solutions to the local interpolation problems, but also on the method used to determine which points lie in a given subdomain.
We use a tree-like data structure (namely, \citeauthor{MATLAB2020a}'s \texttt{KDTreeSearcher} object), as recommended by \cite{Wendland2006}, so the cost of searching the datapoints is $\bigO(\log N)$, with $\bigO(N\log N)$ setup.

\begin{figure}[h]
	\centering
	\includegraphics[width=0.6\linewidth]{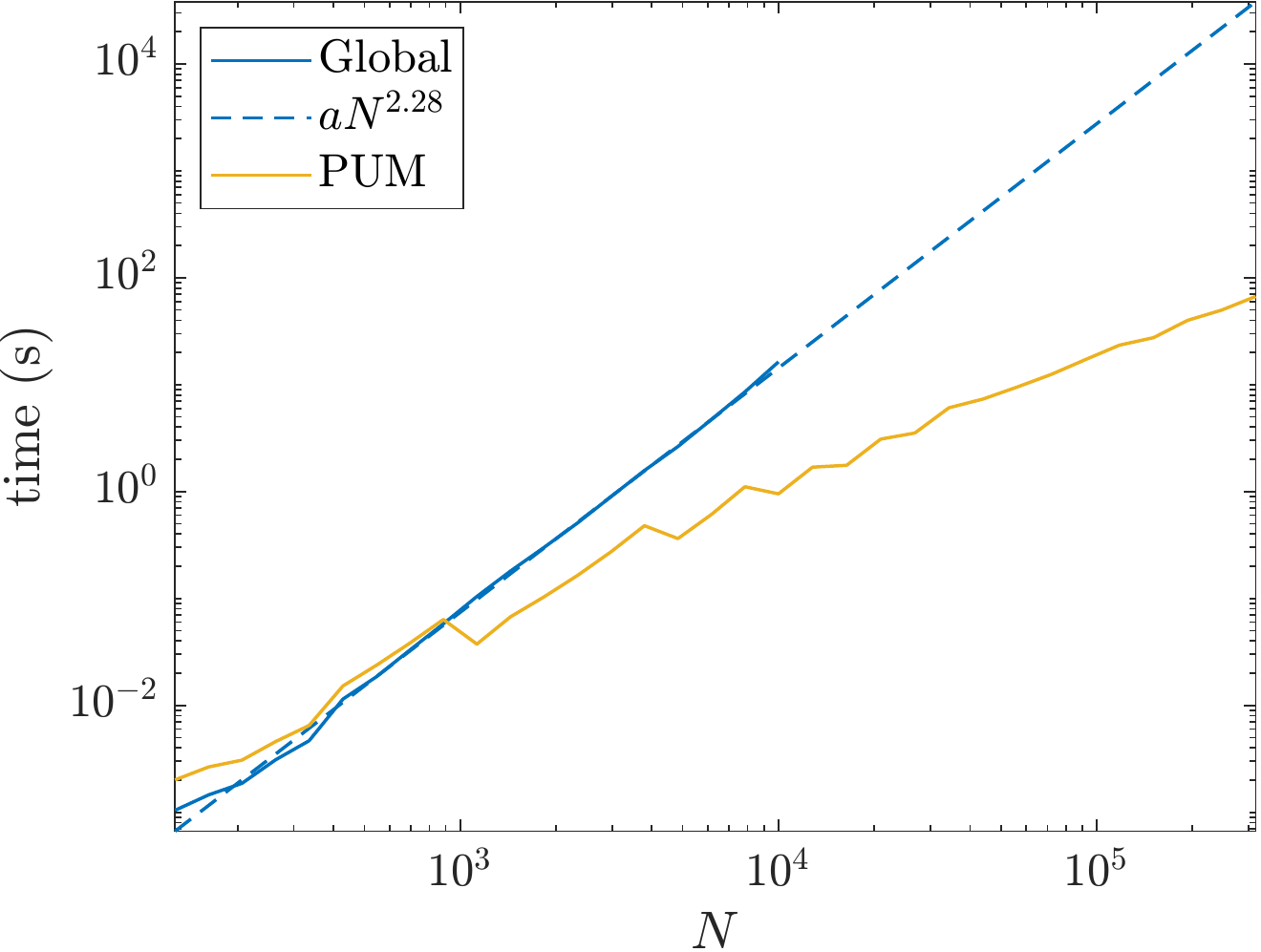}
	\caption[Plot of time taken for global and PU interpolation methods.]{Wall-clock time taken to fit a 2D thin plate spline interpolant to random data (average of 10 trials). For $N>10^4$, the distance matrix for the global interpolant is prohibitively large, so a power law is fit and extrapolated for comparison. Subdomains are chosen to contain approximately $1000$ datapoints each.}
	\label{pumComplexity}
\end{figure}

Figure \ref{pumComplexity} demonstrates the computational advantages of the partition of unity method with around $1000$ datapoints per subdomain. 
The runtime of the global method is not shown for $N>10^4$, as the distance matrix is dense and prohibitively large. 
Note, however, the full distance matrix is never required for the partition of unity method.

\section{Postprocessing} \label{s:post}

\subsection{Domain for evaluation}
\begin{figure}[h]
    \centering
    \begin{subfigure}{0.3\linewidth}
        \includegraphics[width=\linewidth]{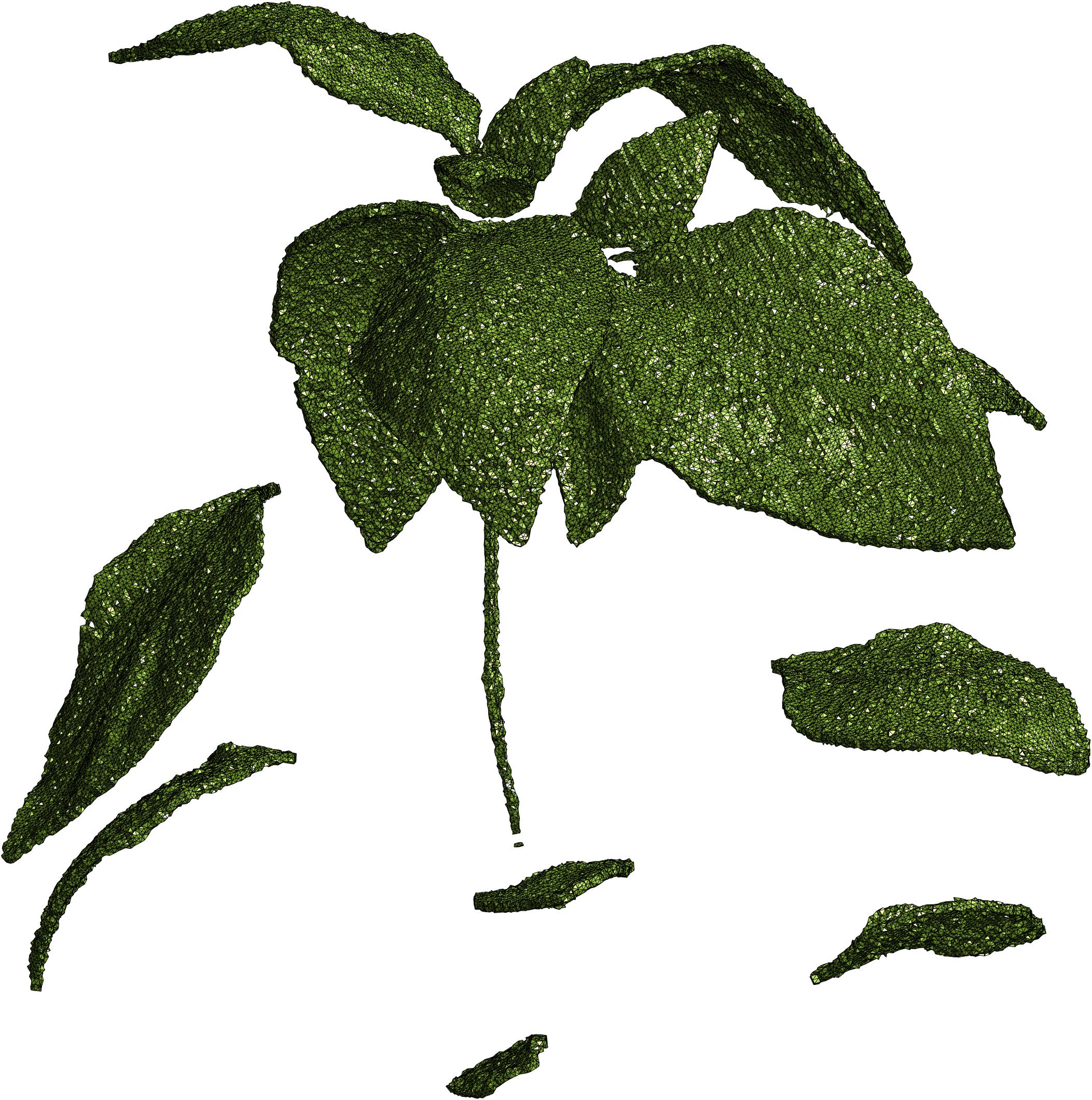}
        \caption{$\alpha = L$}
    \end{subfigure}
    \hspace{10pt}
    \begin{subfigure}{0.3\linewidth}
        \includegraphics[width=\linewidth]{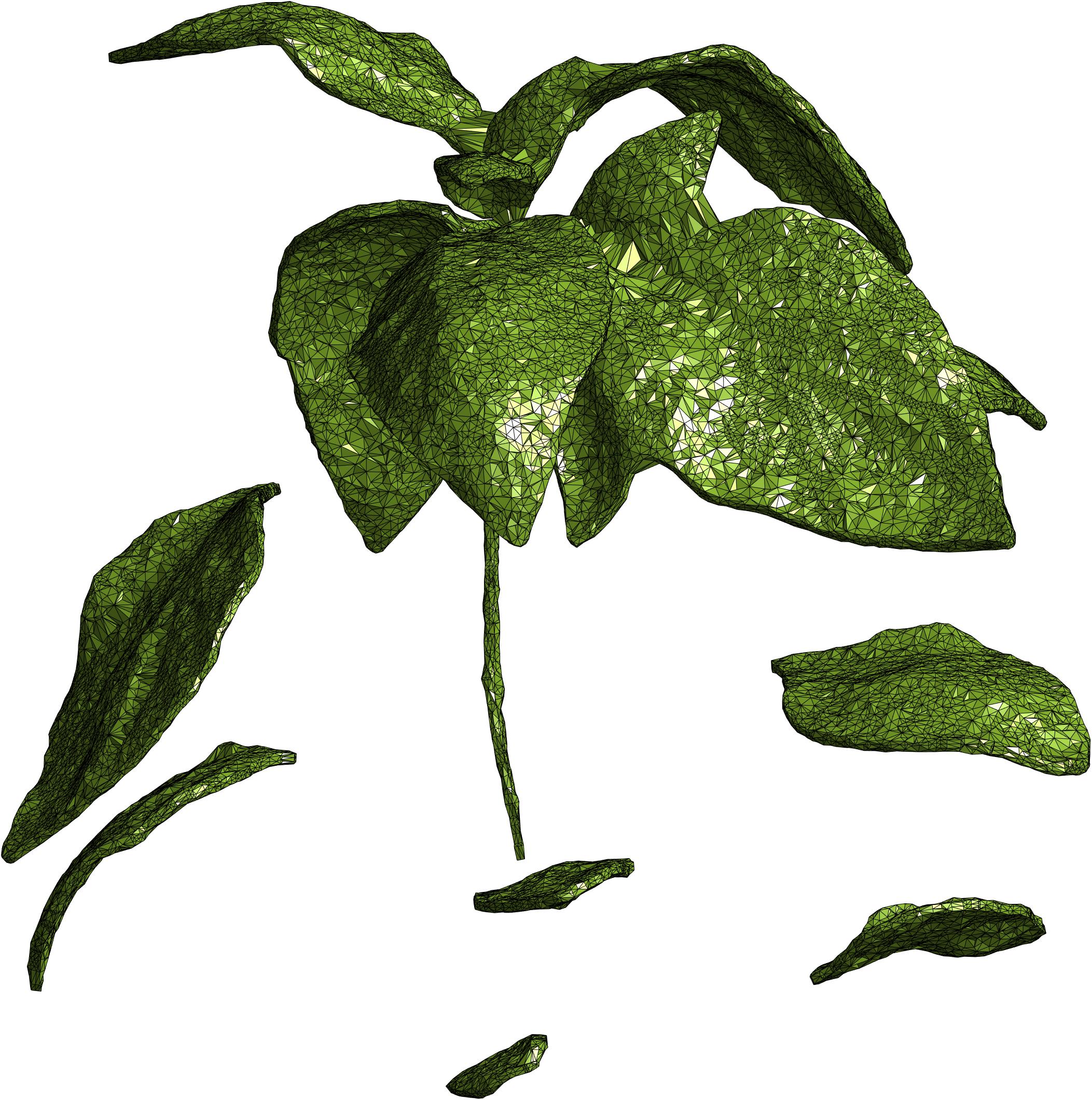}
        \caption{$\alpha = 5L$}
    \end{subfigure}
    \\[10pt]
    \begin{subfigure}{0.3\linewidth}
        \includegraphics[width=\linewidth]{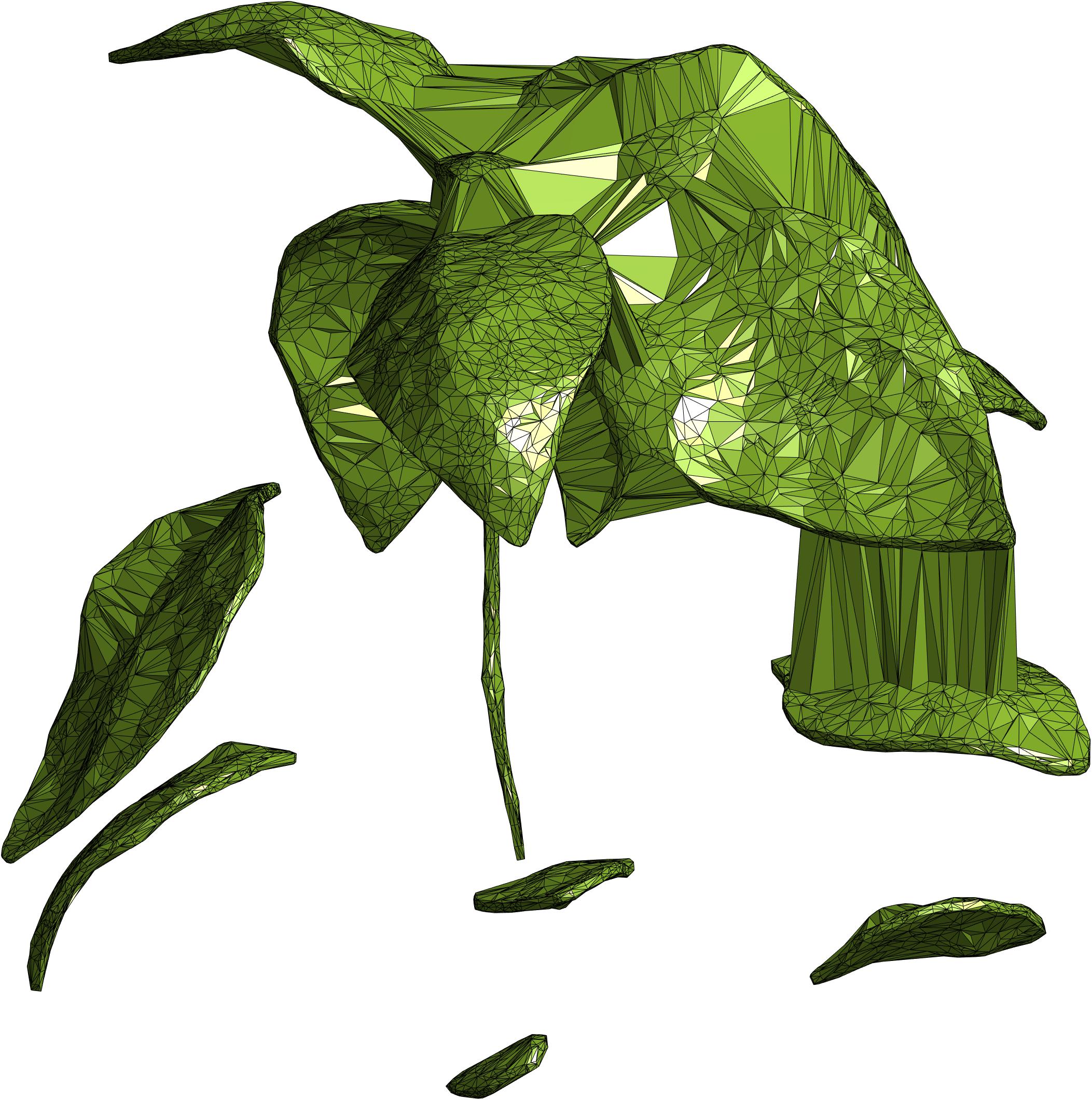}
        \caption{$\alpha = 20L$}
    \end{subfigure}
    \hspace{10pt}
    \begin{subfigure}{0.3\linewidth}
        \includegraphics[width=\linewidth]{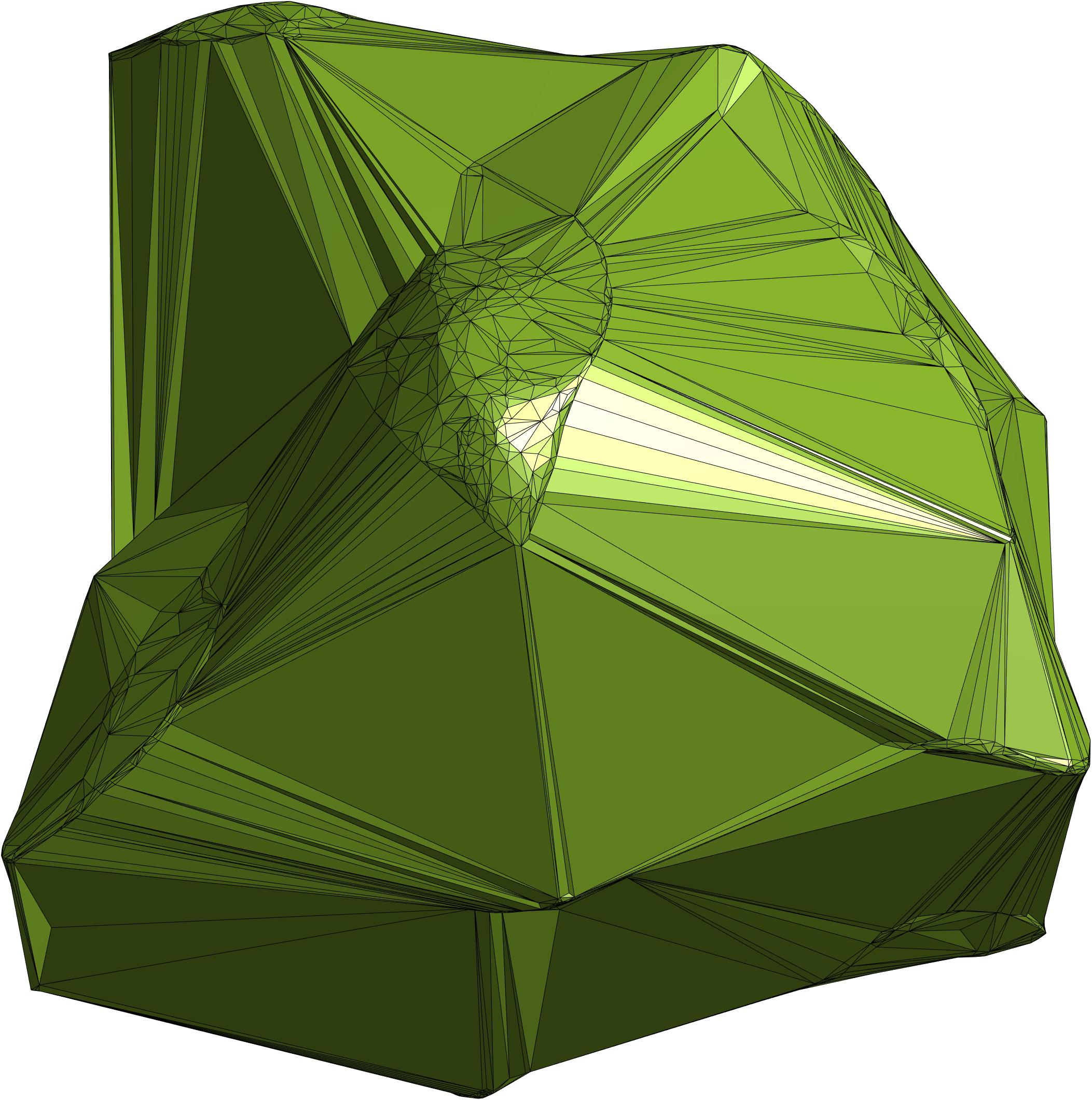}
        \caption{$\alpha = 100L$}
    \end{subfigure}
    \caption{Some $\alpha$-shapes for the capsicum point cloud, associated with radii $\alpha$ that are multiples of $L$, the offset length of the off-surface points. 
    The $\alpha$-shape is a triangulated surface forming a closed volume around the extended point cloud (including off-surface points).}
    \label{fig:alphashape}
\end{figure}

The partition of unity method restricts the definition of the interpolant to the union of the subdomains.
However, we must still take care to restrict the sampling of the interpolant to near the point cloud, as sampling too far from the true surface may result in erroneous surface sheets \citep{Kazhdan2006}. 
We do this using \citeauthor{MATLAB2020a}'s \texttt{alphaShape} \citep{MATLAB2020a}, which, given a radius, $\alpha$, returns a triangulated surface enveloping a set of points \citep{Edelsbrunner1994}. 
A smaller $\alpha$ produces a tighter $\alpha$-shape, while a large $\alpha$ tends to produce a convex hull-like $\alpha$-shape.
Figure \ref{fig:alphashape} shows some $\alpha$-shapes for the capsicum plant point cloud -- note that the triangulations around distinct leaves may not be connected, but are still component-wise watertight.

\begin{figure}
    \centering
    \includegraphics[width=0.8\textwidth]{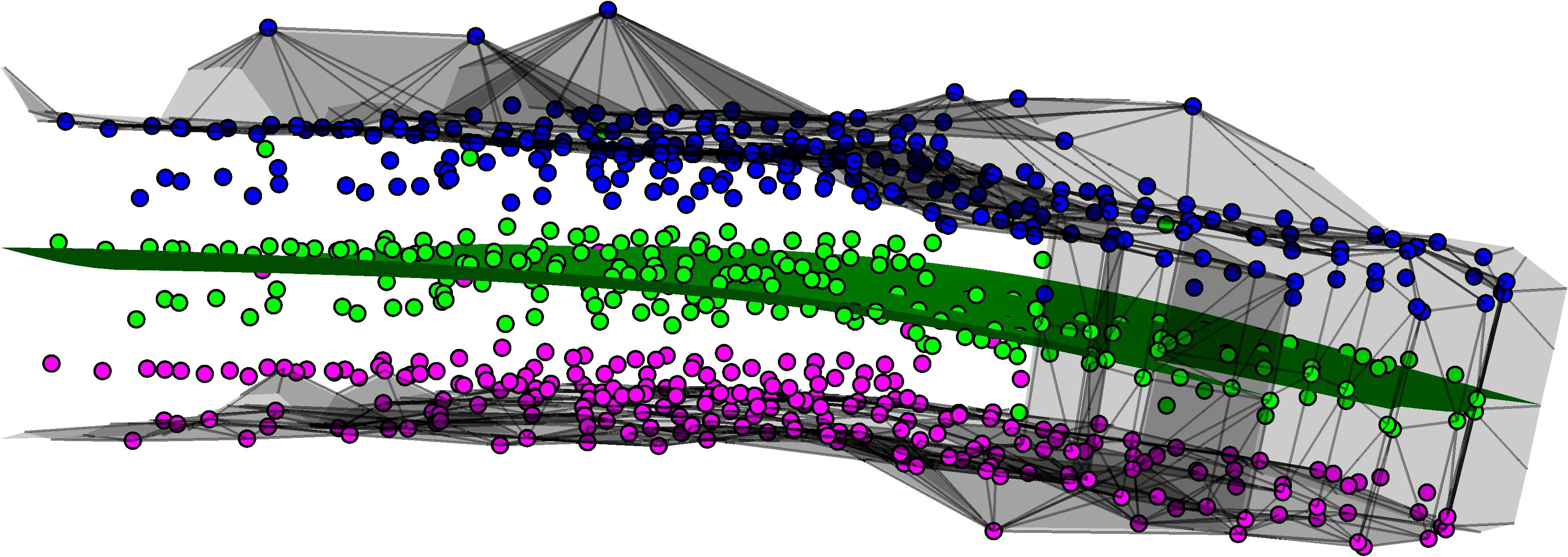}
    \caption{\textcolor{revisedText}{A portion of the domain of a reconstructed leaf. Points show the on-surface (green) and off-surface (magenta and blue) data. The grey mesh is the boundary of the $\alpha$-shape. The green surface is the reconstructed leaf.}}
    \label{fig:pointsSurfAndAshape}
\end{figure}

\textcolor{revisedText}{
Figure \ref{fig:pointsSurfAndAshape} shows a section of the domain with the leaf surface, the on- and off-surface points, and the $\alpha$-shape plotted. 
The $\alpha$-shape forms a tight volume around the points -- the interpolant $\mathcal{F}$ is only defined inside it.
Note that the $\alpha$-shape alone would fail to reconstruct the leaf surface in the presence of noisy data, but is ideally suited to the role of restricting the domain of the implicit function $\mathcal{F}$ to a neighbourhood around the point cloud.
}

We have found that appropriate values for $\alpha$ lie between $3L$ and $10L$, where $L$ is the height of the off-surface points, $\vec{x}_j \pm L\vec{n}_j$.
MATLAB's \texttt{alphaShape} class contains a method to determine whether a point lies outside the $\alpha$-shape.
When sampling the interpolant, we simply exclude such points.

\subsection{Extracting the isosurface}
\begin{figure}[h]
	\centering
	\begin{subfigure}{0.3\textwidth}
	    \includegraphics[width=\linewidth]{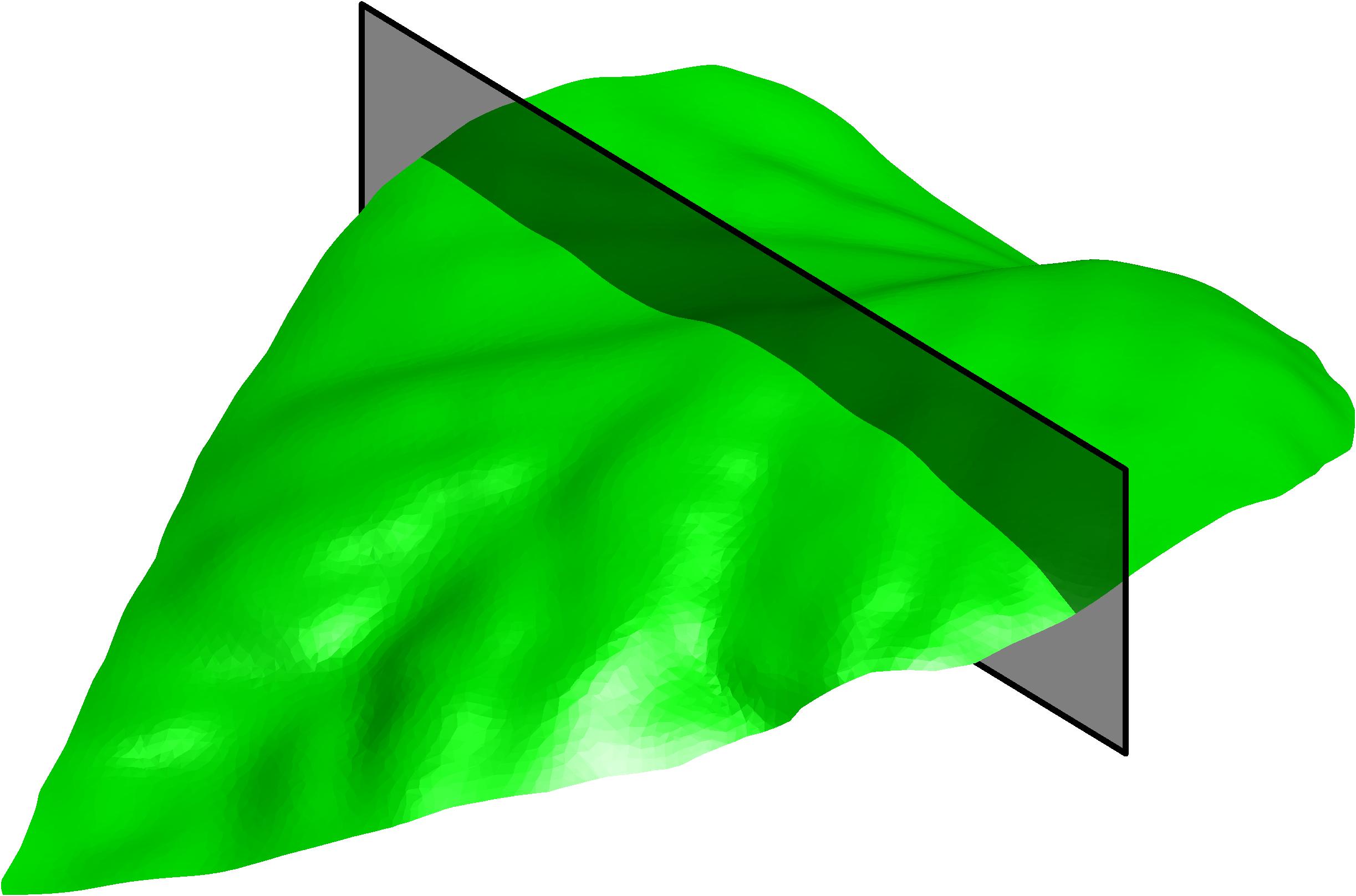}
	\end{subfigure}
	\hspace{10pt}
	\begin{subfigure}{0.6\textwidth}
	    \includegraphics[width=\linewidth]{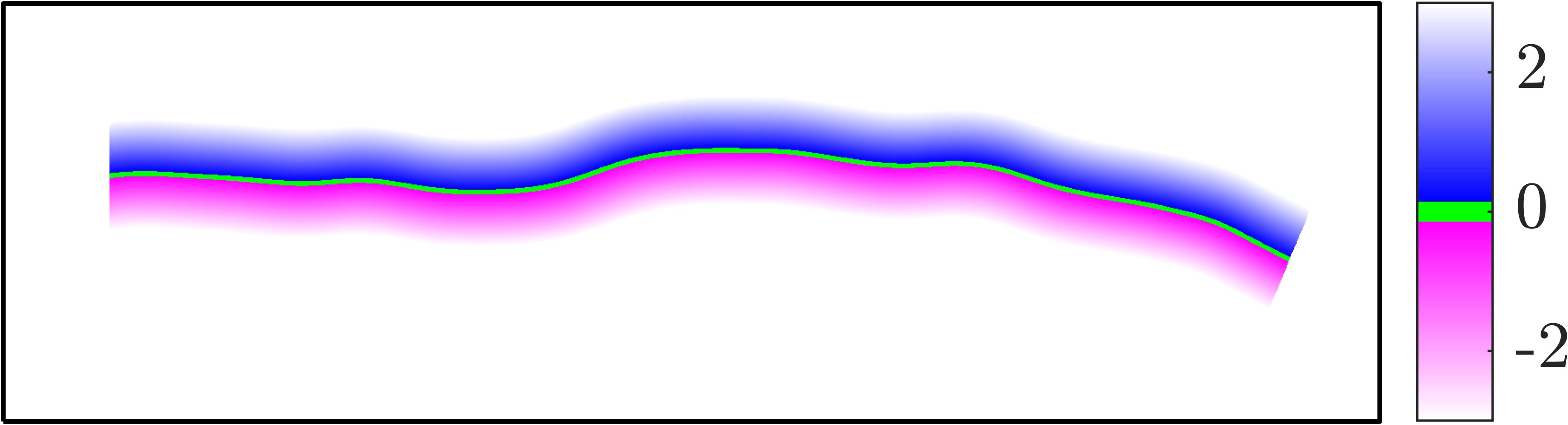}
	\end{subfigure}
	\caption{Left: a reconstruction of a capsicum leaf and a plane through it. Right: the value of the function $\mathcal{F}$ on the plane. The green line where $\mathcal{F}$ is zero valued is the reconstructed leaf surface.}
	\label{fig:levelset}
\end{figure}

If we sample the interpolant on a regular grid, we can extract an approximate surface as a triangulation using \citeauthor{MATLAB2020a}'s \texttt{isosurface} \citep{MATLAB2020a}.
We have also found some success in meshing the $\alpha$-shape with tetrahedra, sampling the interpolant at the vertices, then extracting a triangulated surface with the marching tetrahedra algorithm \citep{doi1991efficient}.
Figure \ref{fig:levelset} shows a triangulated leaf surface, extracted from an implicit RBF-PUM interpolant. 
A slice through the leaf shows the value of the interpolant -- negative on one side of the leaf, zero at the surface, then positive on the other side.
This is the essence of the implicit reconstruction.
\section{Results} \label{s:results}
We have implemented the RBF-PUM in \citeauthor{MATLAB2020a}, and reconstructed plant surfaces using the $C^2$ spline $\phi(r) = r^3$. 
In general, the method is applicable to point clouds with millions of points.
Here, however, we have used grid average downsampled versions of our scanned point clouds, in part to lessen the effect of noise. 
The other advantage of this approach is to improve the conditioning of the interpolation problem: the condition number of the matrix $\vec{A}$ in equation (\ref{smoothAlambdab}) increases exponentially with the minimum seperation distance between datasites \citep{fasshauer2007meshfree}.
By using a grid average downsampling, the minimum seperation distance increases with the grid spacing.

The local interpolation problems are formulated as in section \ref{s:fit} and solved using \citeauthor{MATLAB2020a}'s \texttt{backslash}. 
We avoided the use of a basis-function-specific fast solver here, to demonstrate the flexibility of the RBF-PUM interpolation.
That is, the ability of our method to deal with large datasets is due to the partition of unity method,  which is independent of the local interpolation method.
Furthermore, we note that our implementation has not taken full advantage of parallel computing, although in theory the local interpolation problems are entirely decoupled and may be solved simultaneously.

\begin{figure}
    \centering
    \begin{subfigure}[t]{0.4\textwidth}
        \includegraphics[width=\textwidth]{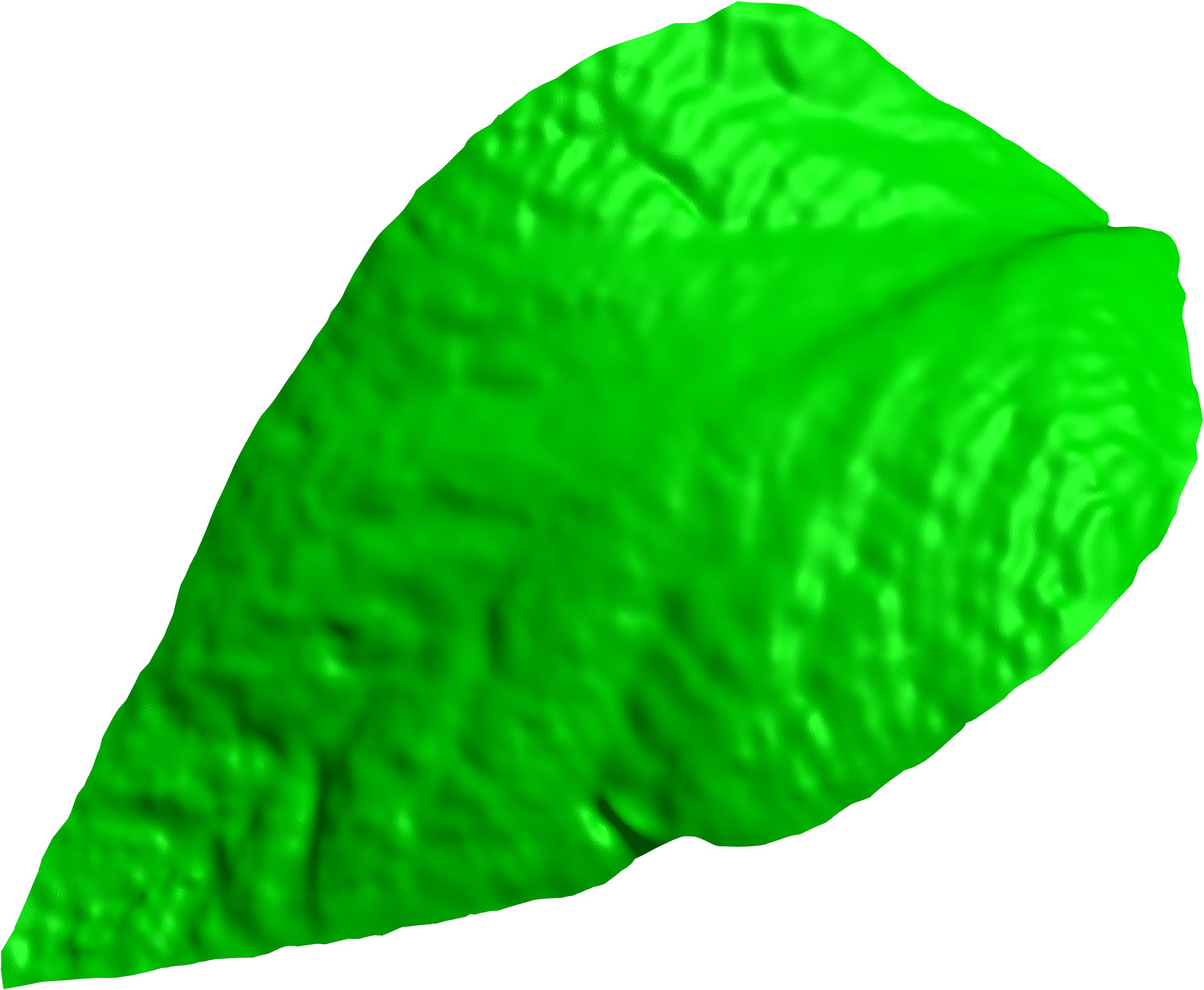}
        \caption{$\rho = 10^{-5}$}
    \end{subfigure}
    \hspace{10pt}
    \begin{subfigure}[t]{0.4\textwidth}
        \includegraphics[width=\textwidth]{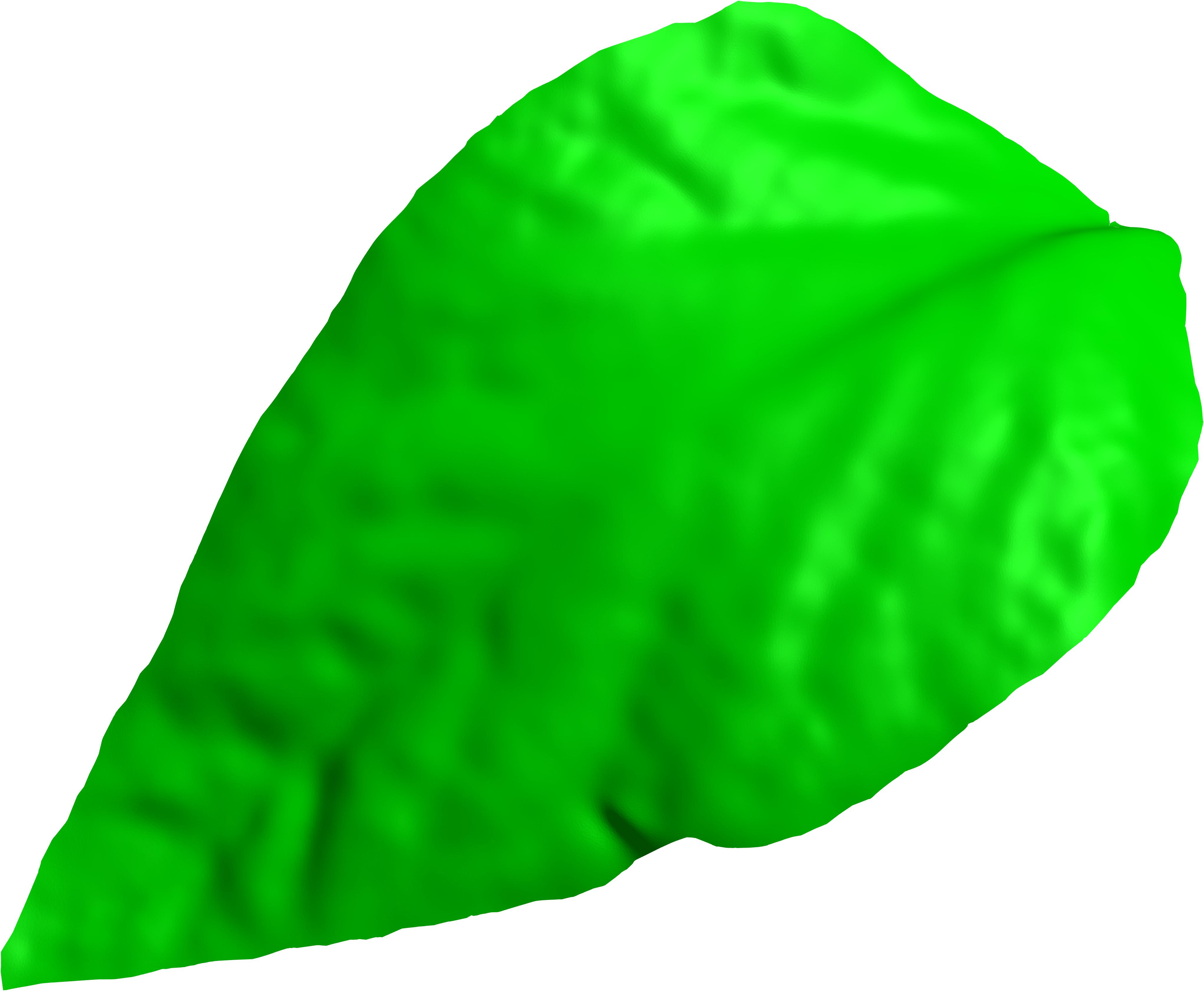}
        \caption{$\rho = 10^{-4}$}
    \end{subfigure}
    \\[10pt]
    \begin{subfigure}[t]{0.4\textwidth}
        \includegraphics[width=\textwidth]{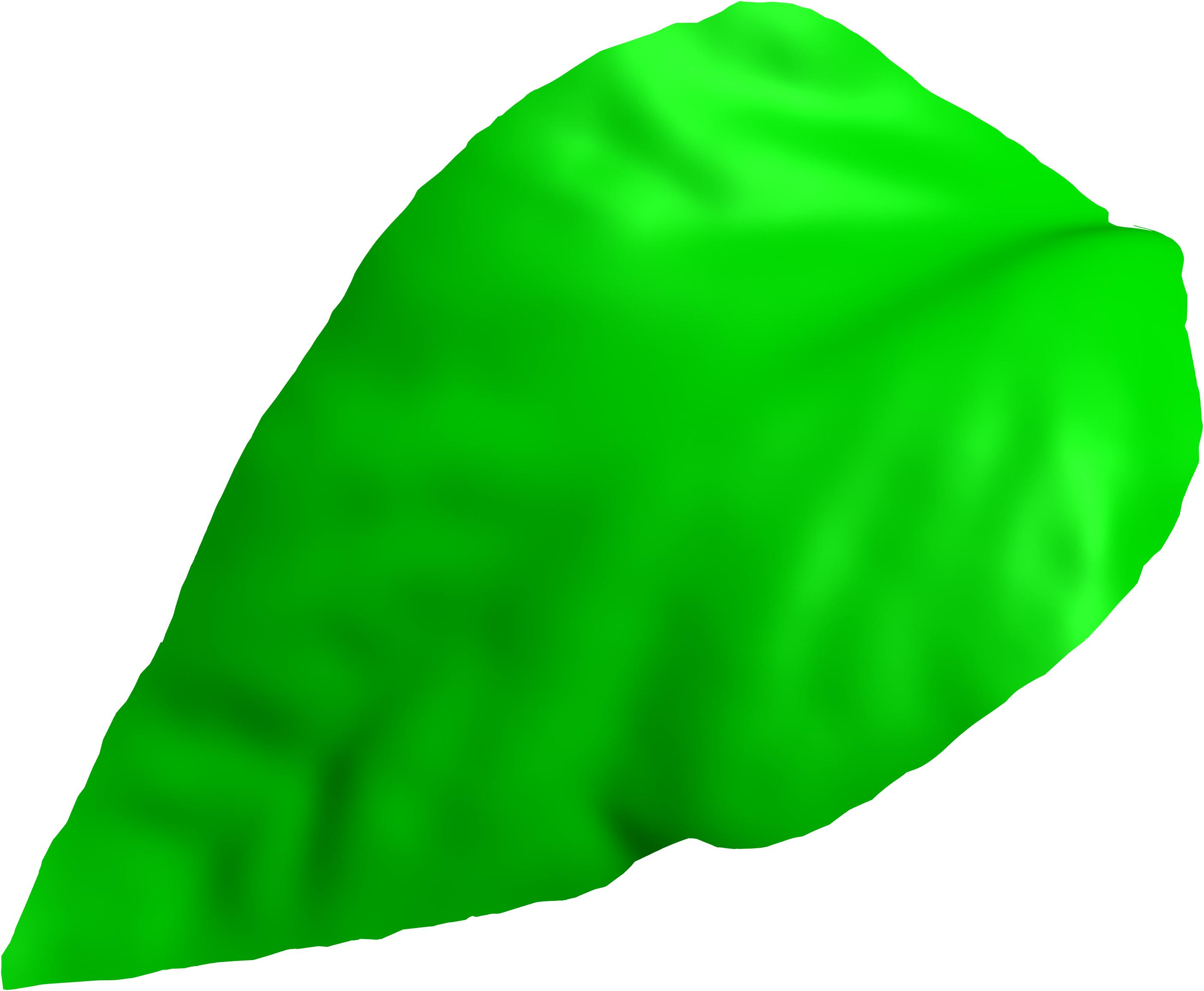}
        \caption{$\rho = 10^{-3}$}
    \end{subfigure}
    \hspace{10pt}
    \begin{subfigure}[t]{0.4\textwidth}
        \includegraphics[width=\textwidth]{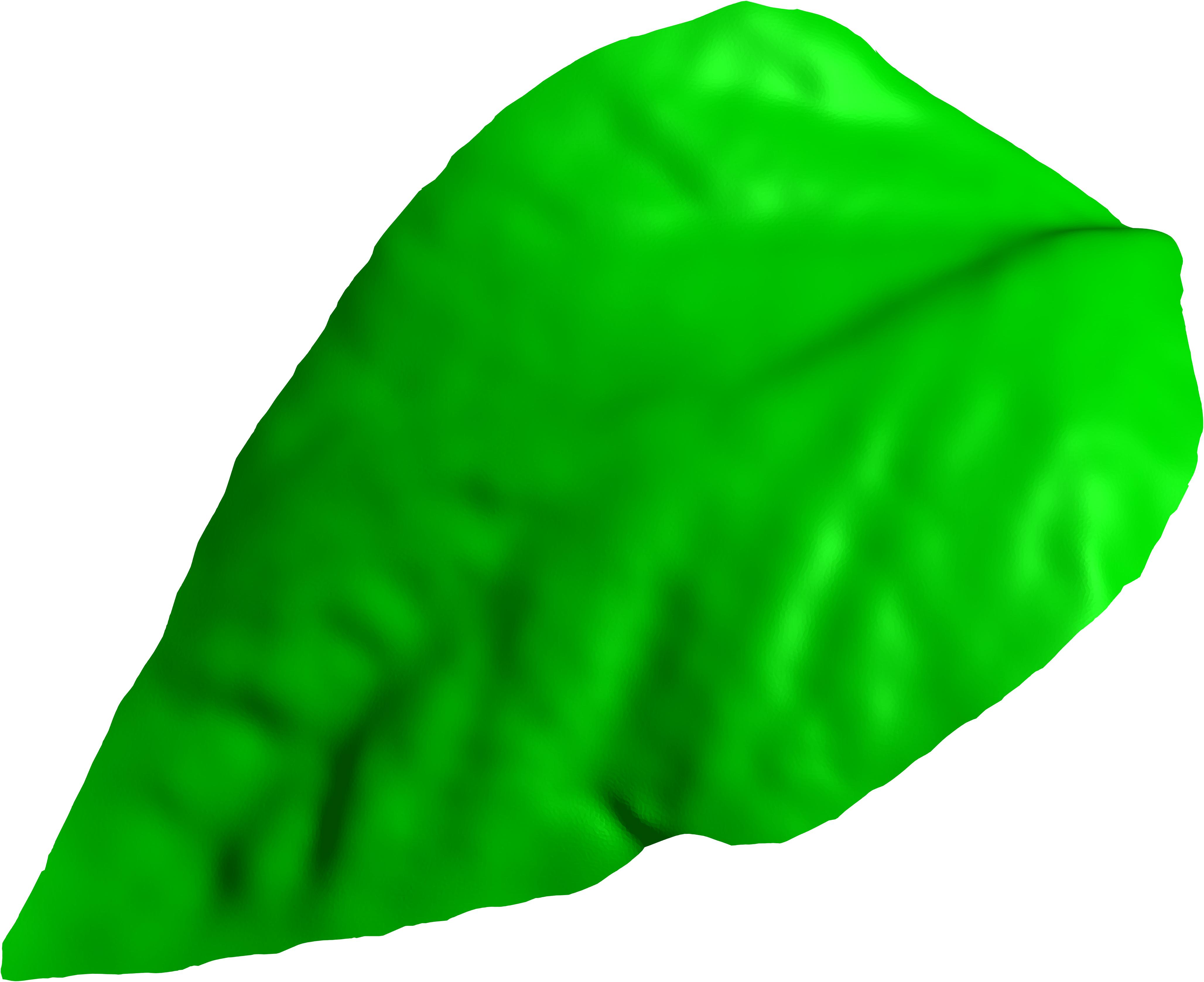}
        \caption{Optimal $\rho$ estimated by GCV for each of 36 subdomains. Median $\rho$ is $3.1153\times10^{-4}$.}
    \end{subfigure}
    \\[10pt]
    \begin{subfigure}[c]{0.4\textwidth}
        \centering
        \includegraphics[width=0.78\textwidth]{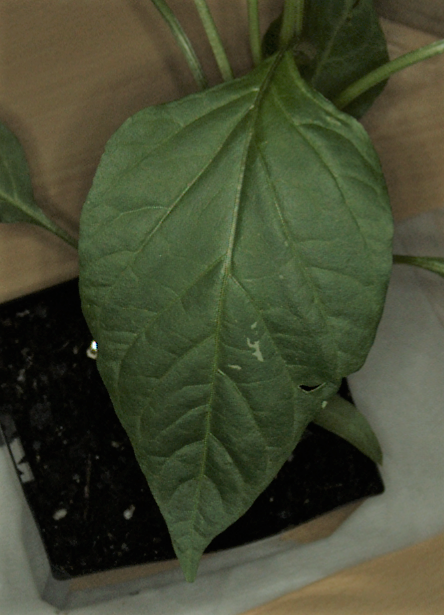}
        \caption{}
    \end{subfigure}
    \hspace{10pt}
    \begin{subfigure}[c]{0.4\textwidth}
        \includegraphics[width=\textwidth]{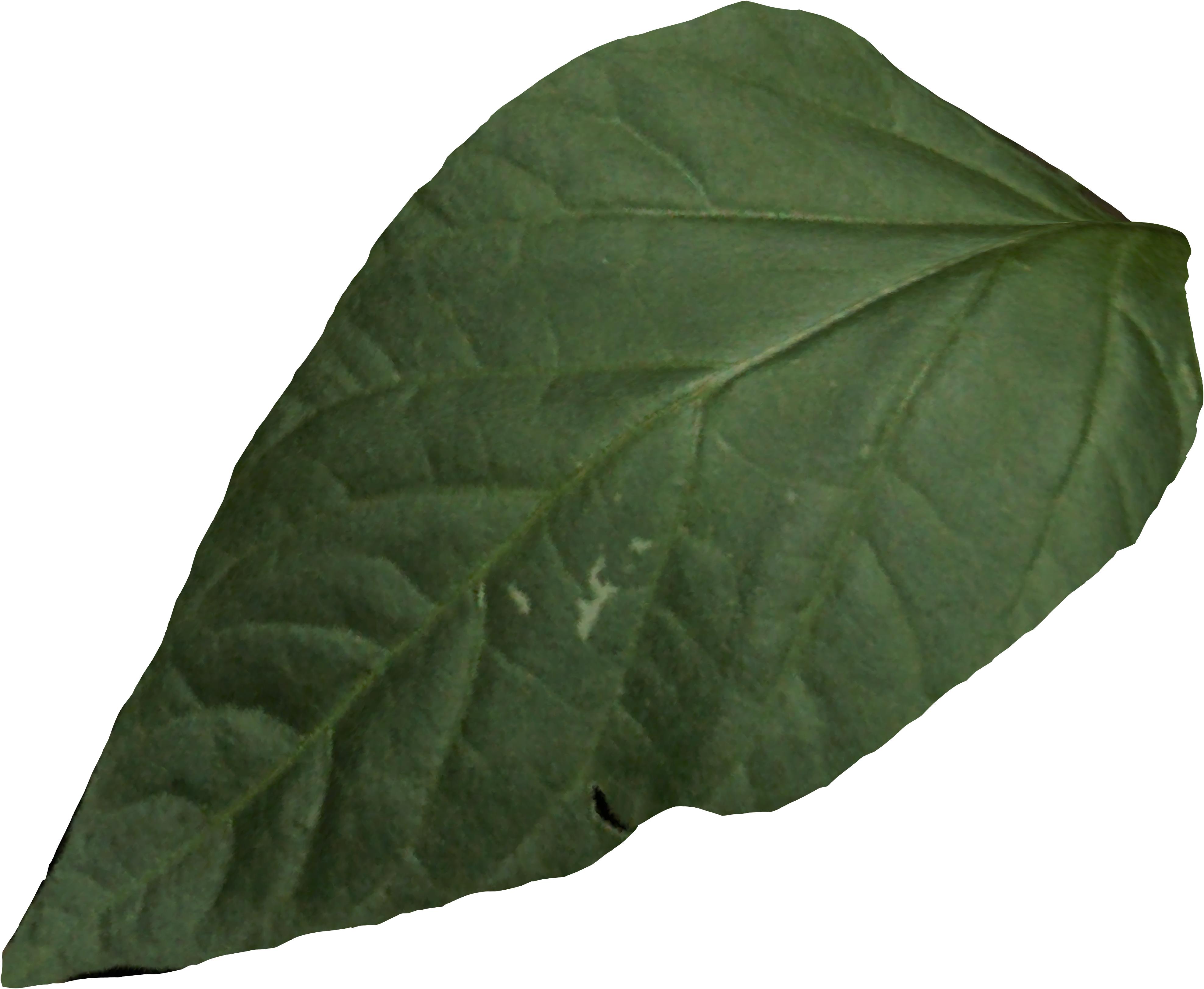}
        \caption{}
    \end{subfigure}
    \caption{Reconstruction of a capsicum leaf. (a-d) A comparison of some reconstructions using global choices of $\rho$ with a reconstruction using locally GCV-estimated $\rho$. (e) A photograph of the capsicum leaf, enhanced for clarity. (f) The GCV-estimate reconstruction, coloured by projecting the image of the leaf onto the surface.}
    \label{fig:rho_comparison}
\end{figure}

By default, we choose a constant value for the smoothing parameter $\rho$ across all of the subdomains.
We believe this choice is suitable given that the noise present is introduced in the scanning process, and should not vary greatly between leaves.
For other applications of this method to datasets with spatially varying noise, generalised cross-validation may be necessary to choose the smoothing parameter on a per-subdomain basis.
Figure \ref{fig:rho_comparison} shows a comparison of reconstructions of a capsicum leaf: three with constant, global smoothing parameters ($\rho$), and one with locally GCV-estimated smoothing parameters. 
We see that the GCV estimate preserves major features of the leaf, but avoids overfitting the noisy scan data.

\begin{figure}[h]
    \centering
    \includegraphics[width=\linewidth]{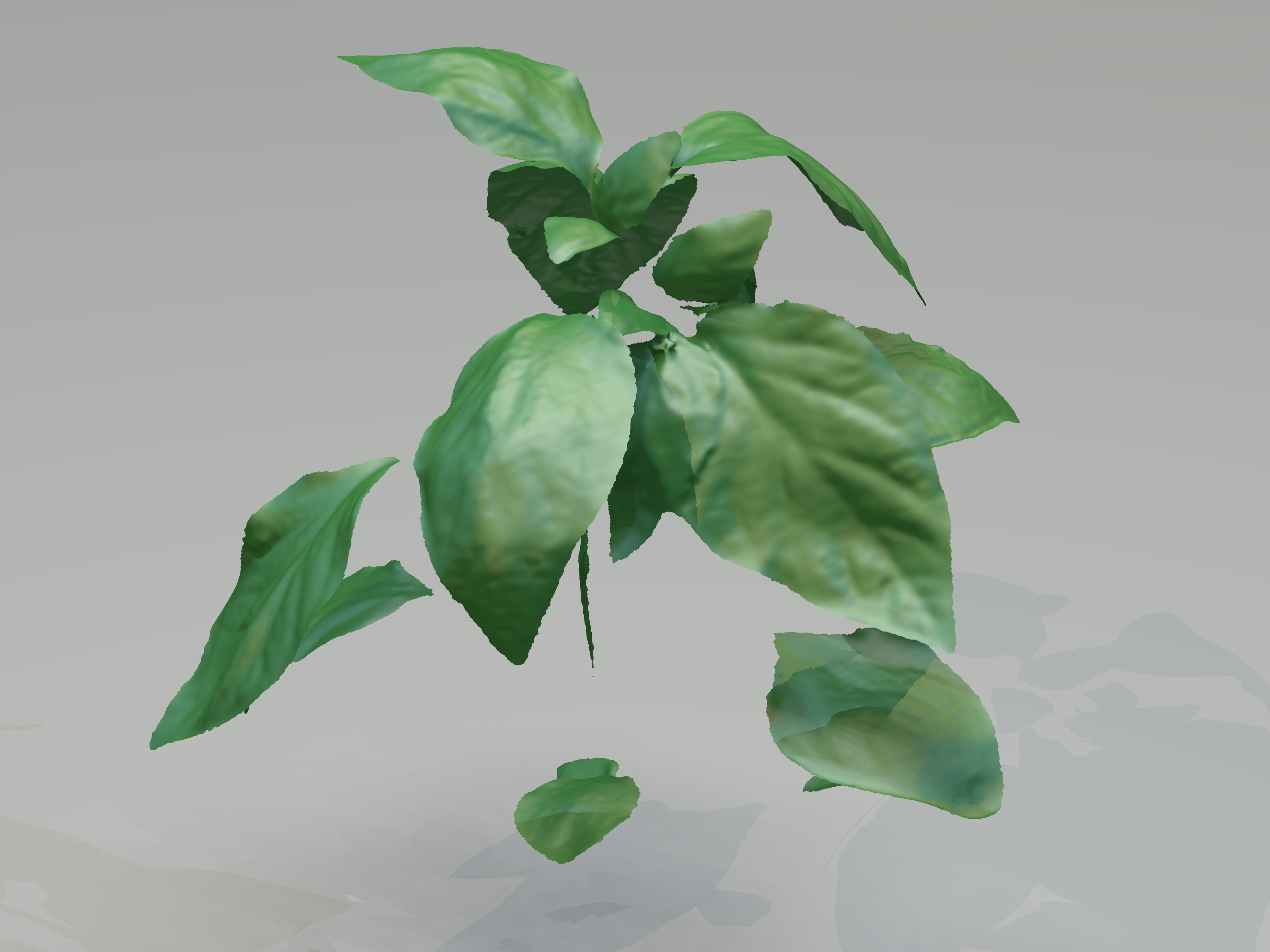}
    \caption{\color{revisedText}A reconstruction of a capsicum plant from a point cloud with $605,112$ points, using $5717$ subdomains. The interpolant was fit in around 1.6 hours, and evaluated on a regular grid (with grid step size 0.25) in around 1.8 hours.
    This image was rendered from the triangulation using the ray-tracing software POV-Ray \citep{povray}.}
    \label{fig:cap}
\end{figure}

\textcolor{revisedText}{Figure \ref{fig:cap} shows a reconstruction of a scanned capsicum plant, rendered using the ray-tracing software POV-Ray \citep{povray}. }
Some artifacts are present near where small leaves protrude from a shared stem, although this is not unexpected: the normal direction changes quite sharply there -- too sharply to be captured in the resolution of the point cloud.
The scan resolution is also the reason that the stem of the plant is only partially captured (where it is thickest).
Our reconstruction captures the vein structure and even boundaries of the larger, separated leaves. 
We note that the boundary is a by-product of the $\alpha$-shape restricting evaluations, and as such is not smooth, despite appearances.
\textcolor{revisedText}{A table of parameter values used in this reconstruction can be found in Table \ref{tab:parms}.}

\subsection{Curvature}
One of the key advantages of our method is the flexibility in the choice of basis function. 
We sought this flexibility so as to use polyharmonic splines with appropriately high orders of continuity, for use in related work on droplet simulation on leaf surfaces.
Specifically, we require continuous mean curvature across our reconstructed surfaces.
We can calculate the mean curvature, $K_M$, of a surface defined implicitly by $\mathcal{F}=0$ \citep{Goldman2005}:
\begin{align*}
    K_M &= -\nabla \cdot \left( \frac{\nabla\mathcal{F}}{\|\nabla\mathcal{F}\|} \right),
\end{align*}
where $\nabla\mathcal{F}$ is the gradient of $\mathcal{F}$.

\begin{figure}[h]
    \centering
    \begin{subfigure}{0.45\linewidth}
    \centering
    \includegraphics[width=\linewidth]{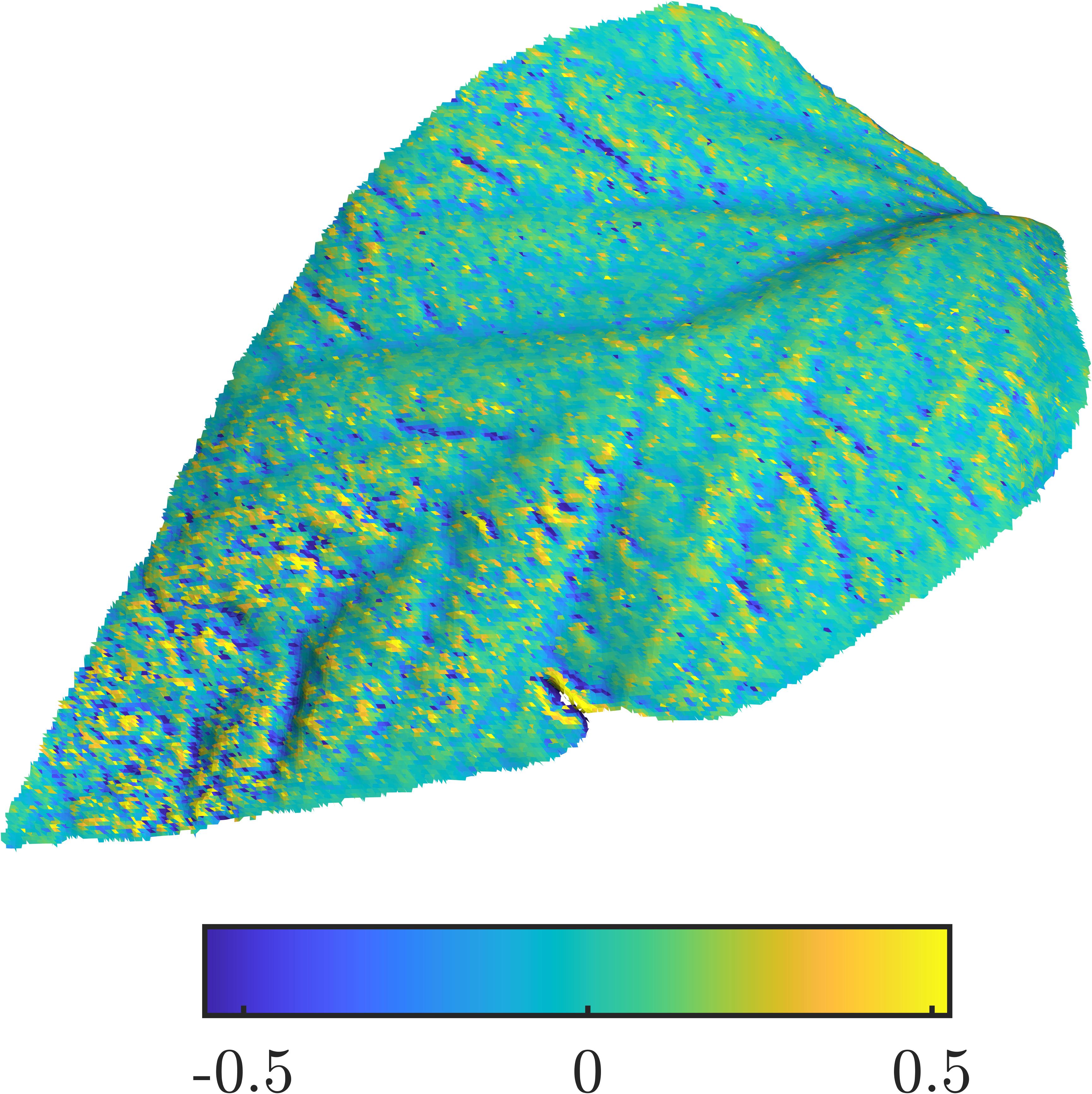}
    \end{subfigure}
    \hfill
    \begin{subfigure}{0.45\linewidth}
    \centering
    \includegraphics[width=\linewidth]{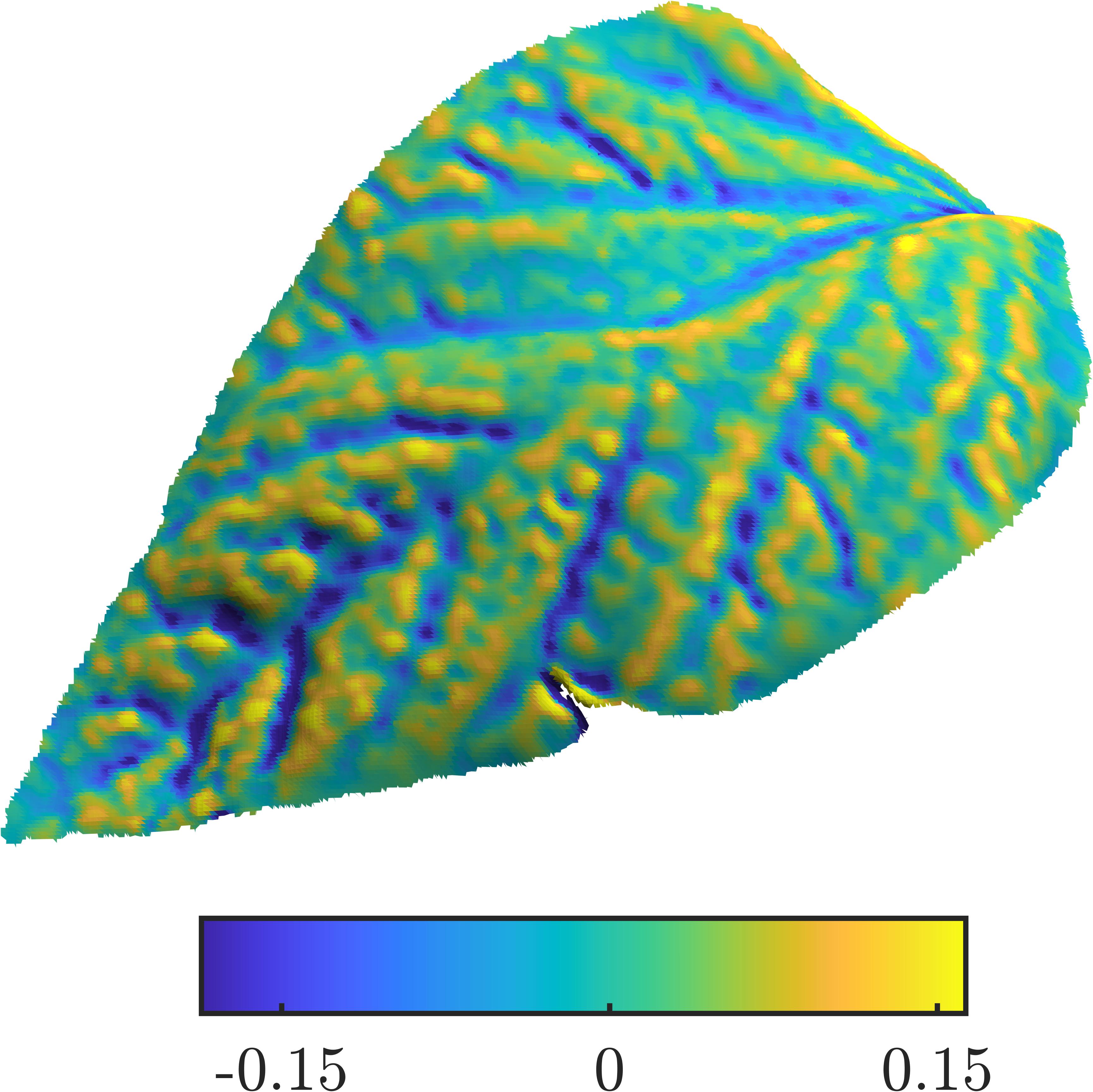}
    \end{subfigure}
    \caption{A reconstruction of a capsicum leaf from a downsampled point cloud of 65,609 points, using a partition of unity with 437 subdomains. Colour indicates mean curvature.
    Left: reconstructed from the $C^0$ spline $\phi(r)=r$, showing discontinuity in mean curvature. 
    Right: reconstructed from the $C^2$ spline $\phi(r)=r^3$, which has continuous second derivatives and thus produces a surface with continuous curvature.}
    \label{fig:curvature}
\end{figure}

Figure \ref{fig:curvature} shows the curvature of a capsicum leaf, reconstructed from the same point cloud using two different polyharmonic splines: $\phi(r) = r$ (which is $C^0$), and $\phi(r) = r^3$ (which is $C^2$).
The discontinuity of second derivatives manifests in the $C^0$ surface as sharp changes in curvature.
The $C^2$ surface, however, does not exhibit these discontinuities.
The importance of this continuous curvature is in future work with these surfaces: for example, thin film models for droplet spreading on leaves will involve a surface curvature term, so we require at least a $C^2$ surface.

\begin{figure}
    \centering
    \includegraphics[width=0.5\linewidth]{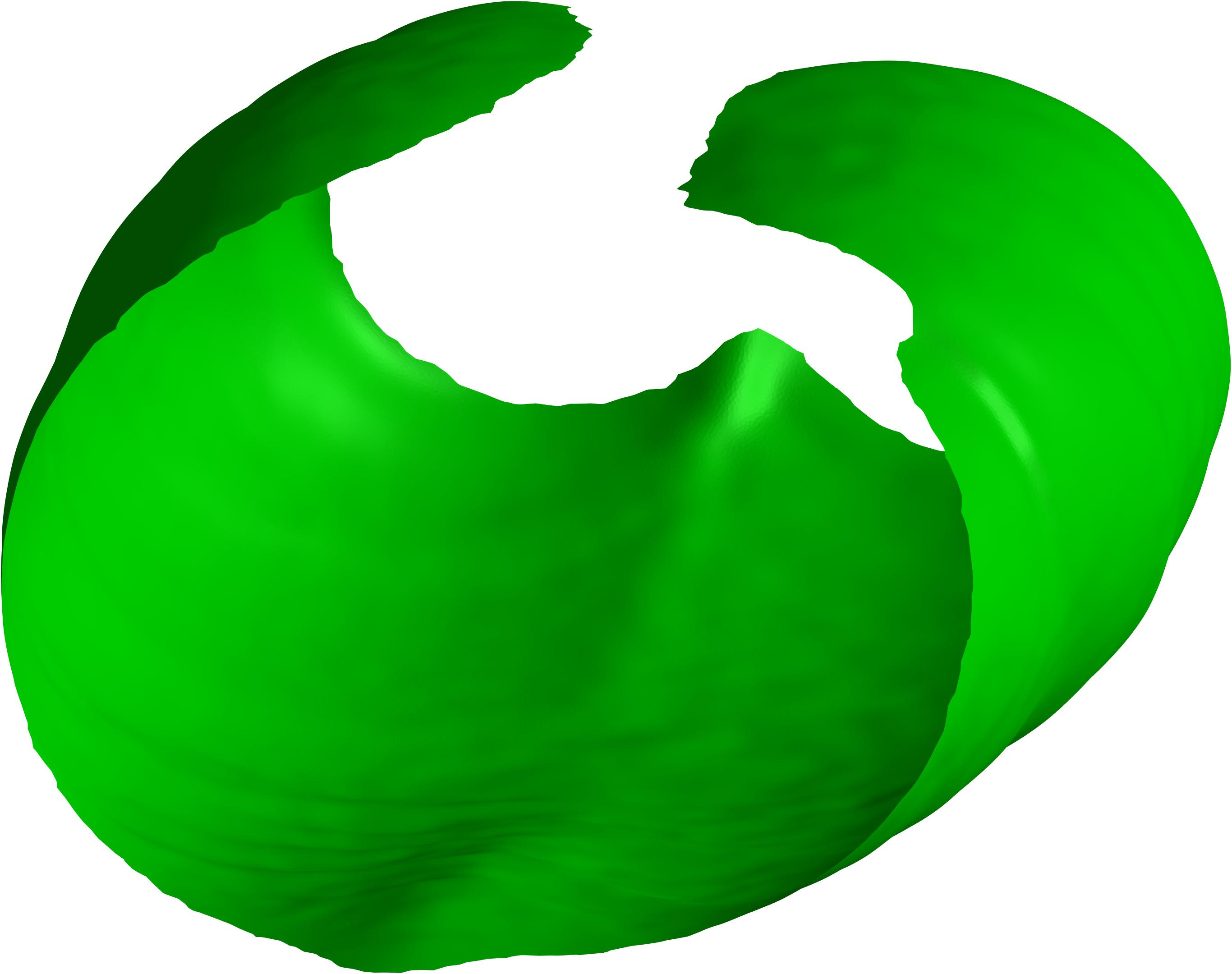}
    \caption{\textcolor{revisedText}{A reconstruction of a transformed point set, demonstrating the method's ability to reconstruct a surface an explicit method could not.}}
    \label{fig:curlyLeaf}
\end{figure}

\textcolor{revisedText}{
\subsection{On curling leaves}
To demonstrate the advantages of our implicit method, we have reconstructed an artificially `curled' leaf. 
We start with a point set $(x_j,y_j,z_j)$ with $-1 < y < 1$, rotated such that the $x$, $y$, and $z$ directions are the three principal components of the point set (in order). 
We then apply the following coordinate transform:
\begin{gather*}
    \beta = \pi \left( y - \frac{1}{2} \right)
    \\
    \eta = \frac{1}{3} \cos ( 2 \beta ) + 1 - z, \\
    x^* = \frac{3x}{2}, \quad
    y^* = \eta \cos(\beta), \quad
    z^* = \eta \sin(\beta).
\end{gather*}
We then reconstruct the leaf surface in the new coordinates, $(x^*_j,y^*_j,z^*_j)$, with no modifications to the method necessary.
Figure \ref{fig:curlyLeaf} shows the reconstruction.
We note that the small gap between the upper edges of the leaf has not been incorrectly bridged, due to the $\alpha$-shape restricting the domain tightly around the point cloud.
}

\section{Conclusion} \label{s:conclusion}
We have described and tested a novel method for reconstructing thin surfaces, built up from existing techniques.
The only condition we impose on the input point cloud is that the true surface have continuously varying normals.
Our interpolation approach combines smoothing polyharmonic splines with a partition of unity method, and an adaptive octree-like process for partitioning the domain.
For evaluation, we use alpha-shapes to define a tight domain around the point cloud to avoid evaluating the interpolant far away from the true surface, where it may be incorrectly zero-valued.

The modular nature of the method means that it would be straightforward to adapt to different local interpolation methods: no fast solver is required.
Furthermore, the decoupling of the local interpolation problems means this method could take advantage of the efficiencies of parallel computing, with only minor modifications.

Going forward, the plant surfaces reconstructed by this method will be incorporated into forthcoming models for droplet spreading and movement on leaves. 
They will also be added to the growing suite of plant models available in a `plant retention spray model' (PRSm) \citep{Dorr2014,Dorr2016,Zabkiewicz2020}.
There they will be used to simulate interception, shatter, bounce and retention of droplets on plant crops with an outlook towards developing more efficient agrichemical spraying processes.
\section*{Acknowledgements}
This research was supported by an Australian Government Research Training Program Scholarship.
The authors gratefully acknowledge support from the Australian Research Council through the ARC Linkage Project (LP160100707) and the associated industry partners Syngenta and Nufarm. 
The capsicum and tomato plants were scanned using the Artex Eva scanner provided by QUT’s eResearch Department. 
Thanks also to Arvind Kumar and Justin Cooper-White for their support in acquiring the capsicum and tomato plants and scanning them at the Australian Institute for Bioengineering and Nanotechnology, University of Queensland.
\textcolor{revisedText}{The authors also thank the anonymous reviewers for their constructive feedback, which improved the final version of the manuscript.} 
\bibliography{thinLeafRefs}
\appendix
\newpage \onecolumn
\section{Pseudocode} \label{pseudocode}

\begin{algorithm}
	\caption{Orienting normals \citep{Hoppe1992}}
	\label{alg:orientNormals}
	\begin{spacing}{1.1}
	\begin{algorithmic}
		\REQUIRE $\mathcal{X}=\{\vec{x}_1,\dots,\vec{x}_N\} \subset \mathbb{R}^s$, $\mathcal{N} = \{\vec{n}_1,\dots,\vec{n}_N\}$, \texttt{coarseGrid}, \texttt{graphNbrs}, \texttt{pcaNbrs}
		\ENSURE $\mathcal{N}$
		\STATE $\mathcal{X}^* \leftarrow \texttt{pcdownsample}(\mathcal{X},\texttt{coarseGrid})$
		\COMMENT \texttt{MATLAB} routine
	    \STATE $\mathcal{N^*} \leftarrow \texttt{pcnormals}(\mathcal{X}^*,\texttt{pcaNbrs})$
	    \COMMENT \texttt{MATLAB} routine
	    \STATE initialise sparse graph $G$ with vertices $\mathcal{X}$
	    \FORALL{ $\vec{x}^*_i$ in $\mathcal{X}^*$ }
	        \STATE find \texttt{graphNbrs} nearest neighbours of $\vec{x}^*_i$ in $\mathcal{X}^*$
	        \STATE exclude neighbours more than $2*\texttt{coarseGrid}$ away
	        \FORALL{ neighbours $\vec{x}^*_j$ }
	            \STATE $G[i,j] \leftarrow 1 - |\vec{n}^*_i \cdot \vec{n}^*_j |$
	        \ENDFOR
	        \STATE $G[i,i] \leftarrow 0$
	    \ENDFOR
	    \STATE $T \leftarrow \texttt{minspantree}(G)$ 
	    \COMMENT \texttt{MATLAB} routine; `sparse' and `forest' options
	    \STATE make list of \texttt{edges} for breadth-first search of $T$ 
	    \COMMENT in \texttt{MATLAB}: \texttt{bfsearch}
	    \FOR{ $k = 1$ \TO number of edges }
	        \STATE $(i,j) \leftarrow \texttt{edges}[k]$
	        \IF{$ \vec{n}^*_i \cdot \vec{n}^*_j < 0 $}
	            \STATE $\vec{n}^*_j \leftarrow -\vec{n}^*_j$ 
	            \COMMENT if normals disagree, flip normal
            \ENDIF
	    \ENDFOR
	    \FORALL{ $\vec{x}^*_i$ in $\mathcal{X}^*$ }
	        \STATE find \texttt{neighbours} of $\vec{x}^*_i$ in $\mathcal{X}$, within radius \texttt{coarseGrid}
	        \FORALL{ $\vec{x}_j$ in \texttt{neighbours} }
	            \IF{$ \vec{n}^*_i \cdot \vec{n}_j < 0 $}
	            \STATE $\vec{n}_j \leftarrow -\vec{n}_j$
            \ENDIF
	        \ENDFOR
	    \ENDFOR
	\end{algorithmic}
	\end{spacing}
\end{algorithm}

\newpage

\begin{algorithm}
	\caption{Octree method for partition of unity}
	\label{PUMoctree}
	\begin{spacing}{1.1}
	\begin{algorithmic}
		\REQUIRE $\mathcal{X}=\{\vec{x}_1,\dots,\vec{x}_N\} \subset \mathbb{R}^s$, $n_\text{min}$, $n_\text{max}$.
		\ENSURE $\mathcal{P}$, $\mathcal{R}$
		\STATE $\vec{y}_1 \leftarrow 
		0.5 \left( 
		[\max_j x_{j,1}, \dots, \max_j x_{j,s}]^T -
		[\min_j x_{j,1}, \dots, \min_j x_{j,s}]^T
		\right)$ 
		\COMMENT{initial cube centre}
		\STATE $\mathcal{P} \leftarrow \{ \vec{y}_1 \}$
		\COMMENT{set of cube centres}
		\STATE $l_1 \leftarrow 2\max_j \norm{\vec{y}_1 - \vec{x}_j}{\infty}$
		\COMMENT{side length of initial cube}
		\STATE $\mathcal{L}\leftarrow\{ l_1 \}$
		\COMMENT{set of side lengths}
		\STATE $r_1 \leftarrow \sqrt{ s\left(\frac{l_1}{2}\right)^2 }$
		\COMMENT{radius of initial covering sphere}
		\STATE $\mathcal{R} \leftarrow \{ r_1 \}$
		\COMMENT{set of radii}
		\WHILE{$\max_i \texttt{countData}(\mathcal{X},\vec{y}_i,r_i) > n_\text{max}$}
			\STATE $i^* \leftarrow \argmax_i \texttt{countData}(\mathcal{X},\vec{y}_i,r_i)$
			\STATE $\mathcal{P}^* \leftarrow \displaystyle\prod_{k=1}^s \left\{ y_{i^*,k} - \frac{r_{i^*}}{2}, y_{i^*,k} + \frac{r_{i^*}}{2} \right\} $
			\COMMENT{set product}
			\STATE $\mathcal{L}^* \leftarrow \left\{ \underbracket[0.5pt]{ \frac{l_{i^*}}{2}, \dots, \frac{l_{i^*}}{2} }_{2^s \text{ times}} \right\}$
			\STATE $\mathcal{R}^* \leftarrow \left\{ \underbracket[0.5pt]{ \sqrt{ s\left(\frac{l_{i^*}}{2}\right)^2 }, \dots, \sqrt{ s\left(\frac{l_{i^*}}{2}\right)^2 }}_{2^s \text{ times}} \right\}$
			\STATE $\mathcal{P} \leftarrow ( \mathcal{P} \setminus \{\vec{y}_{i^*}\} ) \cup \mathcal{P}^*$
			\STATE $\mathcal{L} \leftarrow ( \mathcal{L} \setminus \{l_{i^*}\} ) \cup \mathcal{L}^*$
			\STATE $\mathcal{R} \leftarrow ( \mathcal{R} \setminus \{r_{i^*}\} ) \cup \mathcal{R}^*$
		\ENDWHILE
		\FOR{$i=1$ \TO $|\mathcal{P}|$}
			\STATE $n \leftarrow \texttt{countData}(\mathcal{X},\vec{y}_i,r_i)$
			\IF{$n=0$}
				\STATE mark subdomain $i$ for deletion
			\ELSIF{$n < n_\text{min}$}
				\STATE find $\vec{x}^*$, the $(n_\text{min})^\text{th}$ nearest neighbour of $\vec{y}_i$ in $\mathcal{X}$. 
				\STATE $r_i \leftarrow \norm{\vec{y}_i - \vec{x}^*}{2}$
			\ENDIF
		\ENDFOR
		\STATE delete subdomains marked for deletion
	\end{algorithmic}
	\end{spacing}
\end{algorithm}

\section{\textcolor{revisedText}{Capsicum plant reconstruction parameters}}
\setcounter{table}{0}

\begin{table}[h]
    \renewcommand{\arraystretch}{1.2}
    \color{revisedText}
    \centering
    \caption{\color{revisedText}The parameters used to reconstruct the capsicum plant shown in Figure \ref{fig:cap}.}
    \label{tab:parms}
    \begin{tabular}{ccp{0.6\linewidth}}
        \toprule
        Parameter & Value & Description \\ \midrule
        \texttt{denoiseNbrs} & 50 & number of neighbours used to determine if a point is an outlier to be discarded 
        \\
        \texttt{denoiseThreshold} & 0.15 & how many standard deviations must a point's average distance to neighbours be from the mean to be classified an outlier
        \\
        $\Delta x$ & 0.5 & grid step for grid-average downsampling
        \\
        \texttt{pcaNbrs} & 50 & number of neighbours used to approximate surface normals
        \\
        $\Delta x_\text{coarse}$ & 2.0 & coarser downsampling grid step, only for the purposes of forming the normal similarity graph
        \\
        \texttt{graphNbrs} & 10 & number of neighbours used to build weighted graph of normal similarities
        \\
        $L$ & $2\Delta x$ & off-surface point distance
        \\
        $N_\text{min}$ & 2000 & minimum number of points in a subdomain
        \\
        $N_\text{max}$ & 5000 & maximum number of points in a subdomain
        \\
        \texttt{expand} & 1.1 & factor by which to multiply the initial radius of a new subdomain to ensure some overlap between neighbouring subdomains
        \\
        $\phi(r)$ & $r^3$ & polyharmonic spline radial basis function
        \\
        $\rho$ & $10^{-3}$ & smoothing parameter
        \\
        $\alpha$ & $2\Delta x$ & $\alpha$-shape `tightness' parameter
        \\
        \bottomrule
    \end{tabular}
\end{table}

\end{document}